\documentclass[twocolumn]{autart}    

\usepackage{amsmath, amssymb}
\usepackage[authoryear,longnamesfirst]{natbib}
\usepackage{hyperref}
\usepackage{url}
\setcitestyle{authoryear,open={(},close={)}} 
\hypersetup{
    hidelinks,
    colorlinks=true,
    linkcolor=blue,
    citecolor=blue
}
\newtheorem{remark}{Remark}

\newtheorem{ass}{Assumption}

\usepackage{graphicx}          

\begin{document}

\begin{frontmatter}

\title{Fully Coupled Nonlinear FBS$\bigtriangleup$Es: Solvability and LQ Control Insights\thanksref{footnoteinfo}} 

\thanks[footnoteinfo]{This paper was not presented at any IFAC
meeting. Corresponding author Maoning Tang. Tel. +86 572 2321591.
Fax +86 572 2321591.}

\author[Huzhou]{Zhipeng Niu}\ead{zpengniu@163.com},    
\author[Huzhou]{Qingxin Meng}\ead{mqx@zjhu.edu.cn},               
\author[HongKong]{Xun Li}\ead{li.xun@polyu.edu.hk}, 
\author[Huzhou]{Maoning Tang}\ead{tmorning@zjhu.edu.cn}
\address[Huzhou]{Department of Mathematical Sciences, Huzhou University, Zhejiang 313000,PR  China}  
\address[HongKong]{Department of Applied Mathematics, The Hong Kong Polytechnic University, Hong Kong 999077, PR  China}             


\begin{abstract}                          
In this paper, a class of fully coupled nonlinear forward-backward stochastic difference equations (FBS$\bigtriangleup $Es) is proposed and the existence of solutions is proved based on a linear-quadratic (LQ)  optimal control problem. Inspired from the solvability studies of various forward-backward stochastic differential equations (FBSDEs), the dominant-monotone framework is discretised and a continuum approach is used to prove the unique solvability of the fully coupled FBS$\bigtriangleup $Es and to obtain a pair of estimates on the solutions, and finally, the conclusions are applied to the related LQ problem.
\end{abstract}

\end{frontmatter}

\section{Introduction} \label{sec:1}

Since the seminal work of \cite{pardoux1990adapted}, the theory of backward stochastic differential equations (BSDEs) has undergone significant development, owing to their rich mathematical structure and wide range of applications.
When a BSDE is considered together with a forward stochastic differential equation (FSDE), the resulting system is referred to as a forward-backward stochastic differential equation (FBSDE). Such a system is called \emph{coupled} if the forward equation depends on the backward component or vice versa. In particular, if this dependence is mutual---i.e., both the forward and backward equations explicitly depend on each other---the system is called a \emph{fully coupled} FBSDE.

	The time discretization of  FBSDEs  gives rise to forward-backward stochastic difference equations (FBS$\Delta$Es). With the advancement of computational algorithms, FBS$\Delta$Es provide substantial advantages in numerical computation, particularly in the development of efficient solution schemes.
For instance,
\cite{ma2002numberical}  proposed a numerical approximation method of backward stochastic
differential equations;
\cite{zhang2004numerical} designed a numerical scheme for a class of backward stochastic differential equations (BSDEs) with possible path-dependent terminal values;
\cite{gobet2005regression} developed a new numerical scheme based on iterative regressions on function bases, whose coefficients are evaluated using Monte Carlo simulations;
\cite{delarue2006forward} proposed a time-space discretization scheme for quasi-linear parabolic PDEs;
\cite{bender2008time} considered a numerical algorithm to simulate high-dimensional coupled FBSDEs under weak coupling or monotonicity conditions;
\cite{cheridito2013bsdeltaes} provided existence results and comparison principles for solutions of BS$\Delta$Es and then proved convergence of these to solutions of BSDEs;
\cite{gobet2015analytical} analyzed analytical approximations of backward SDEs in the limit of small nonlinearity and short time, in the case of nonsmooth drivers.

    In this paper, we consider a class of FBS$\Delta$Es.
	In fact, fully coupled FBS$\Delta$Es arise naturally in discrete-time stochastic optimal control problems, where the forward equation models the system dynamics and the backward equation represents the cost variable or the value function.
	In finance, such models describe portfolio optimization and asset-liability management, where the forward wealth dynamics and the backward risk-adjusted valuation are inherently interdependent.
	This paper aims to investigate the solvability of a class of fully coupled nonlinear FBS$\Delta$E systems. Before proceeding with the analysis, we provide a motivating example. \cite{khallout2019risk} previously studied the continuous-time version of this example.

	An insurance company manages policyholders' funds with the dual objective of achieving asset growth through investment and meeting the cash flow obligations of its clients. We aim to formulate a discrete-time model to describe this wealth evolution process and analyze the optimal investment strategy.

	Let \( x_k \) denote the wealth of the policyholder at time \( k \), with \( x_0 = m_0 \) being the initial investment. Define the time set as \( \mathbb{T} = \{0, 1, \dots, N-1\} \). At each time \( k \in \mathbb{T} \), the insurance company is required to make a payment of \( c_k x_k \) to the policyholder, where \( c_k \) is the payout ratio.
	The present value of these cash flows, discounted to time \( k \) using a deterministic discount factor \( \exp\left(-\sum_{i=0}^k \lambda_i \right) \), where \( \lambda_i \geq 0 \) is bounded and deterministic, is given by
	\[
	y_k = \mathbb{E}\left[\sum_{j=k}^{N} e^{-\sum_{i=0}^j \lambda_i} c_j x_j \mid \mathcal{F}_k\right].
	\]
	Here, \( y_k \) represents the present value of all future policy payments from time \( k \) to \( N \), a standard concept in actuarial science and financial mathematics. It allows the insurer to evaluate whether current assets are sufficient to cover future liabilities.
	Suppose the insurer invests in a financial market consisting of a risk-free asset with interest rate \( r_k \), and a risky asset that follows a discrete-time geometric Brownian motion with drift \( \mu_k \) and volatility \( \sigma_k \). All coefficients are assumed to be deterministic, bounded, and \( \sigma_k \geq \varepsilon > 0 \). The dynamics of the wealth process are described by:
	\begin{equation} \label{eq:400}
		\left\{
		\begin{aligned}
			x_{k+1} &= (1 + r_k) x_k + \rho_k u_k + \sigma_k u_k \omega_k, \\
			x_0 &= m_0,
		\end{aligned}
		\right.
	\end{equation}
	where \( u_k \) denotes the investment in the risky asset and \( \rho_k = \mu_k - r_k \) is the risk premium.

	However, \( x_k \) only captures the dynamics of forward wealth. From the insurer's perspective, we model the liability valuation using a recursive backward equation:
	\begin{equation} \label{eq:500}
		y_k = \mathbb{E}\left[y_{k+1} (1 + \lambda_k) - c_k x_k \mid \mathcal{F}_k\right],
	\end{equation}
	where \( y_k \) represents the amount the insurer needs to prepare at time \( k \) to meet future obligations. Note that \( y_N \) can be specified as 0 or some terminal liability function.
	Combining \eqref{eq:400} and \eqref{eq:500}, we obtain the following forward-backward stochastic difference system:
	\begin{equation} \label{eq:600}
		\left\{
		\begin{aligned}
			x_{k+1} &= (1 + r_k) x_k + \rho_k u_k + \sigma_k u_k \omega_k, \\
			y_k &= \mathbb{E}\left[(1 + \lambda_k) y_{k+1} - c_k x_k \mid \mathcal{F}_k\right], \\
			x_0 &= m_0, \\
			y_N &= 0.
		\end{aligned}
		\right.
\end{equation}
	Clearly, system \eqref{eq:600} is a fully coupled FBS$\Delta$E, where the backward component \( y_k \) depends on \( x_k \). Such coupling reflects the interdependence between investment performance and liability obligations.

	In this paper, we aim to investigate the solvability of the following nonlinear fully coupled FBS$\Delta$Es:
\begin{equation} \label{eq:1.1}
	\left\{\begin{aligned}
		& x_{k+1}  =b\left(k, \theta_k\right)+\sigma(k, \theta_k) \omega_k, \\
		& y_k = -f\left(k+1, \theta_k\right),\quad k \in \mathbb{T}, \\
		& x_0  = \Lambda (y_0), \\
		& y_N = \Phi(x_N),
	\end{aligned}\right.
\end{equation}
where $\theta_k:=(x_k,y'_{k+1},z'_{k+1}) $, $z_{k+1}:=y_{k+1}w_k, y'_{k+1}:=\mathbb E[y_{k+1}|\mathcal{F}_{k-1}], z'_{k+1}:=\mathbb E[z_{k+1}|\mathcal{F}_{k-1}] $. Moreover, $x(\cdot)$ and $y(\cdot)$ take values in some given Euclidean space $\mathbb{H}$.
Also, $b(\cdot)$, $\sigma(\cdot)$ and $f(\cdot)$ are three given mappings with values in $\mathbb{H}$ and $\left\{w_k\right\}_{k\in \mathbb T}$ is a $\mathbb{R}$-value martingale difference sequence defined on a probability space $(\Omega, \mathcal{F}, P)$ . To simplify the notation,  for any $ \theta=(x,y',z'),$
define
{\small
	\begin{equation}
		\Gamma (k,\theta ) :=(f(k+1,\theta) ,b(k,\theta ),\sigma (k,\theta ))
\end{equation}}
and
{\small
	\begin{equation}
		\langle\Gamma(k, \theta), \theta\rangle :=\langle f(k+1, \theta), x\rangle+\langle b(k, \theta), y'\rangle+\langle\sigma(k, \theta), z'\rangle .
	\end{equation}
	In this way, all coefficients of FBS$\bigtriangleup$Es (\ref{eq:1.1}) are collected by $(\Lambda,\Phi,\Gamma)$.
	
	Coupled FBSDEs were first studied by \cite{antonelli1993backward}, who proved existence and uniqueness for small time intervals and showed by counterexample that the same may not hold for longer durations, even under Lipschitz conditions.
		Recently, a growing body of innovative research has tackled the challenge of solving FBSDEs on arbitrary time domains, mainly through two prevailing approaches. The first approach is called a four-step scheme approach, which combines partial differential equations with probabilistic methods and requires that all coefficients in the equations are non-stochastic and the coefficients of the diffusion terms are non-degenerate. \cite{ma1994solving} first proposed the method, which was further developed by \cite{delarue2002existence} and \cite{zhang2006wellposedness}, etc.
	The second method, called a method of continuation, originates from the stochastic Hamiltonian system arising from solving control problems and obtains the solvability of the random coefficients of the FBSDE for an arbitrary length of time by requiring certain monotonicity conditions on the coefficients. The method was originally proposed by \cite{hu1995solution}, and then further developed by \cite{peng1999fully} and \cite{yong1997finding}.

	We note that in most previous studies of FBSDEs, the initial and terminal states are not coupled, i.e., $\Lambda$ is independent of y and $\Phi$ is independent of $x$; see \cite{antonelli1993backward}, \cite{delarue2002existence}, \cite{hu1995solution},
	\cite{ma1994solving}, \cite{ma2015well} and \cite{ma1999forward} for related studies.
	However, the FBSDEs studied by \cite{li2018forward}, \cite{yong2010forward}, and \cite{yu2022forward} are general cases where coupled initial and terminal states exist. This paper is more concerned with a class of fully coupled FBS$\bigtriangleup$Es.

	Some results have been achieved in the study of linearly fully coupled FBS$\bigtriangleup$Es. \cite{xu2017solvability,xu2018general} have studied a class of fully coupled FBS$\bigtriangleup$Es from the point of view of LQ optimal control, whose equations in which the inverted variables have the form of conditional expectations and satisfy that the coefficient matrices of the inverted variables are degenerate. For fully coupled infinite horizon FBS$\bigtriangleup$Es, \cite{xu2017solution} and \cite{song2020forward} adopted a generalized Riccati equation to systematically describe and explicitly express the adaptation solution. It is worth mentioning that \cite{zhang2016complete} advanced the development of discrete-time FBS$\bigtriangleup$Es by deriving the maximum principle for discrete-time mean-field linear quadratic (LQ) control problems for the first time.
	
 Hower, few rigorous solvability results exist specifically for nonlinear FBS$\Delta$Es.
	\cite{ji2024solvability} established pioneering solvability results for fully coupled FBS$\Delta$Es. Their work provides: (1) a complete characterization of linear FBS$\Delta$E solvability through the almost sure invertibility of matrices $\Gamma_t(P_{t+1})$ with explicit solutions; (2) the first existence and uniqueness theorem for nonlinear FBS$\Delta$Es  under classical monotonicity conditions, resolving discrete-time analysis challenges via a novel continuation method adaptation; and (3)
foundational analysis of BS$\Delta$Es where solutions are obtained through the Galtchouk-Kunita-Watanabe decomposition, yielding existence, uniqueness, and stability estimates

 Regarding FBSDEs, inspired by various stochastic linear LQ problems, \cite{yu2022forward} proposed a more general class of domination-monotonicity framework and used the method of continuation to obtain a unique solvability result and a pair of estimates for general coupled FBSDEs. It is worth mentioning that the domination-monotonicity framework established by \cite{yu2022forward}, by introducing various matrices, matrix-valued random variables, and matrix-valued stochastic processes, can be applied to various stochastic LQ problems. This new framework contains the previous traditional monotonicity condition and develops it.

		In this paper, we establish a class of discrete-time fully coupled nonlinear FBS$\Delta$Es based on a linear quadratic (LQ) optimal control framework by employing a discrete method. In particular, we discretize the domination-monotonicity framework proposed in \cite{yu2022forward} and apply the method of continuation to derive a pair of estimates for the FBS$\Delta$E system \eqref{eq:1.1}, thereby proving its unique solvability. To illustrate the applicability of our theoretical results, we also consider two LQ optimal control problems.
		It is important to emphasize several key differences between our work and previous studies:

		\begin{itemize}
			\item[(i)] Compared with \cite{xu2017solvability,xu2018general}, which systematically studied the solvability of linear FBS$\Delta$Es by establishing explicit relations between the forward and backward components via generalized algebraic Riccati equations and deriving necessary and sufficient conditions for existence and uniqueness, our work considers more general nonlinear fully coupled systems.

			\item[(ii)] In contrast to \cite{song2020forward} and \cite{zhang2016complete}, which deal with general linear FBS$\Delta$Es and express the adapted solutions via Riccati equations, our approach addresses \textit{nonlinear} systems. We construct separate estimates for the solutions and establish uniqueness using the continuation method. The development of numerical methods for such nonlinear systems constitutes a promising direction for future research.
			
			\item[(iii)] While \cite{yu2022forward} studied the unique solvability of fully coupled nonlinear FBSDEs under a new domination-monotonicity condition in continuous time, our work focuses on the discrete-time setting. In this framework, classical stochastic analysis tools---such as It\^{o} formula, quadratic variation, and stochastic integrals---are no longer applicable. Instead, we utilize conditional expectations, martingale difference sequences (MDS), and recursive techniques. Specifically, backward dynamics are handled using recursions such as:
			\[
			\mathbb{E}[\omega_k \mid \mathcal{F}_k] = 0, \quad \mathbb{E}[\omega_k^2 \mid \mathcal{F}_k] = 1, \quad y'_k = \mathbb{E}[y_{k+1} \mid \mathcal{F}_{k+1}],
			\]
			which avoid the path regularity assumptions required in continuous-time models.
			
			\item[(iv)] The class of fully coupled nonlinear FBS$\Delta$Es is also investigated in \cite{ji2024solvability}.  Our formulation differs significantly in the system's structure and the proof techniques.  \cite{ji2024solvability} considered a backward equation with three components $(Y, Z, N)$, whereas our model involves only a single backward component $Y$. Moreover, while \cite{ji2024solvability} relied on the classical monotonicity condition. In contrast, our study adopts the generalized domination-monotonicity condition \cite{yu2022forward}, ?guaranteeing solution existence and uniqueness while extending the applicability of the classical condition.
		\end{itemize}
	}

	The remainder of the paper is organized as follows. In Section~\ref{sec:2}, we introduce the necessary notations and present two preliminary lemmas on stochastic difference equations (S$\Delta$Es) and backward stochastic difference equations (BS$\Delta$Es), which will serve as foundational tools for our analysis. Section~\ref{sec:3} is devoted to the study of the fully coupled FBS$\Delta$E system \eqref{eq:1.1} under domination-monotonicity conditions via the continuation method, where we establish the existence and uniqueness of solutions, along with a pair of crucial estimates, as stated in Theorem~\ref{thm:3.1}. In Section~\ref{sec:4}, we apply the theoretical results to two linear-quadratic (LQ) optimal control problems and demonstrate the existence and uniqueness of the corresponding optimal controls.
Finally, we conclude the paper by summarizing the main results and highlighting potential directions for future research.

\section{  Notations and Preliminaries  } \label{sec:2}
Let $N$ be a given positive integer.
	$\mathbb{T}$ , $\overline{\mathbb{T}}$  and  \b{${\mathbb T}$}   denote the sets  $\{0,1, \cdots, N-1\}$,  $\{0,1,2, \ldots, N\}$ and $\{1,2, \ldots, N\}$ respectively.
	Let $(\Omega, \mathcal{F}, \mathbb{F}, \mathbb{P})$ be a complete filtered probability space
	with a filtration $\mathbb{F}=\left\{\mathcal{F}_k : k=0, \cdots, N-1\right\}$. Assume that $\left\{w_k\right\}_{k\in \mathbb T}$ is a $\mathbb{R}$-value martingale difference sequence defined on a probability space $(\Omega, \mathcal{F}, P)$, i.e.
	$$
	\mathbb{E}\left[w_{k+1} \mid \mathcal{F}_k\right]=0, \quad \mathbb{E}\left[w_{k+1}(w_{k+1})^{\top} \mid \mathcal{F}_k\right]=1.
	$$
	In this paper, we assume that $\mathcal{F}_k$ is the $\sigma$-algebra generated by $\left\{x_0, w_l, l=\right.$ $0,1, \ldots, k\}$. For convenience, $\mathcal{F}_{-1}$ denotes $\left \{ \emptyset,\Omega  \right \}$.

	Let $\mathbb{R}^n$ be the $n$-dimensional Euclidean space with the norm $|\cdot|$ and the inner product $\langle\cdot,\cdot\rangle$. Let $\mathbb{S}^n$ be the set of all symmetric matrices in $\mathbb{R}^{n\times n}$. Let $\mathbb{R}^{n\times m}$ be the collection of all $n\times m$ matrices with the norm $|A|=\sqrt{\textrm{tr}(AA^\top)}$, for $\forall A\in \mathbb{R}^{n\times m}$ and the inner product:
	\begin{equation}
		\left\langle A,B\right\rangle = \textrm{tr}(AB^\top),\quad A,B\in \mathbb{R}^{n\times m}.\nonumber
	\end{equation}

	Let $\mathbb H$ be a Hilbert space with norm $\|\cdot\|_\mathbb H$, then we introduce some notations as follows:



	$\bullet$ $L^{2} _{\mathcal{F}_{N-1} }(\Omega ; \mathbb{H})$: the set of all $\mathbb H$-valued  $\mathcal{F}_{N-1}$-measurable random variables $\xi$ satisfying

	\begin{equation}
		\|\xi\|_{L^{2} _{\mathcal{F}_{N-1} }(\Omega ;\mathbb{H})} :=\Big[\mathbb{E}\|\xi \|_\mathbb H^{2} \Big]^{\frac{1}{2}}<\infty .\nonumber
	\end{equation}

	$\bullet$ $L^{\infty} _{\mathcal{F}_{N-1} }(\Omega ;\mathbb H)$: the set of all $\mathbb H$-valued
	$\mathcal{F}_{N-1}$-measurable essentially bounded variables.

	$\bullet$ $L_{\mathbb{F}}^2\left(\mathbb{T}; \mathbb{H}\right)$: the set of all
	$\mathbb{H}$-valued  stochastic process $f(\cdot)=\{f_k|f_k$ is $\mathcal{F}_{k-1}$-measurable, k$\in \mathbb{T}\}$ satisfying

	\begin{equation}
		\|f(\cdot)\|_{L_{\mathbb{F}}^2\left(\mathbb{T}; \mathbb{H}\right)} :=\bigg [\mathbb{E}\bigg (\sum_{k=0}^{N -1} \|f_k\|_{\mathbb{H}}^{2}\bigg ) \bigg]^{\frac{1}{2}}<\infty.\nonumber
	\end{equation}

	$\bullet$ $L_{\mathbb{F}}^{\infty}(\mathbb{T}; \mathbb{H})$: the set of all
	$\mathbb{H}$-valued  essentially bounded stochastic processes $f(\cdot)=\{f_k|f_k$ is $\mathcal{F}_{k-1}$-measurable, k$\in \mathbb{T}\}$.

	$\bullet$ $\mathbb{U}$ : the set of all admissible control
	\begin{align*}
		\mathbb{U} := &\Bigg\{ \mathbf{u} = \Big( u_0, u_1, \ldots, u_{N-1} \Big) \mid
		u_k \text{ is } \mathcal{F}_{k-1} \text{-measurable}, \\
		\quad \quad\quad& u_k \in \mathbb{R}^m, \mathbb{E}\left[ \sum_{k=0}^{N-1} \left| u_k \right|^2 \right] < \infty \Bigg\}.
	\end{align*}

	For the sake of simplicity of notation, we will also present some product space as follows:

	$\bullet$ $N^{2}(\overline{\mathbb{T}};\mathbb{R}^{2n} ):= L_{\mathbb{F}}^2\left(\overline{\mathbb{T}}; \mathbb{R}^n\right)\times L_{\mathbb{F}}^2\left(\overline{\mathbb{T}}; \mathbb{R}^n\right)$ . For any $(x(\cdot ),y(\cdot ) ) \in N^{2}(\overline{\mathbb{T}} ;\mathbb{R}^{2n} )$, its norm is given by
	\begin{equation}
		\begin{aligned}
			\|(x(\cdot),y(\cdot)) \|_{N^{2}(\overline{\mathbb{T}};\mathbb{R}^{2n} )}:= \left \{ \mathbb{E}\bigg [\displaystyle\sum_{k=0}^N|x(k)|^{2}+\displaystyle\sum_{k=0}^N|y(k)|^{2} \bigg ] \right \}^{\frac{1}{2}}.\nonumber
		\end{aligned}
	\end{equation}

	$\bullet$ $\mathcal{N}^{2}(\mathbb{T};\mathbb{R}^{3n} ):= L_{\mathbb{F}}^2\left(\mathbb{T}; \mathbb{R}^n\right)\times L_{\mathbb{F}}^2\left(\mathbb{T}; \mathbb{R}^n\right)\times L_\mathbb{F}^2\left(\mathbb{T}; \mathbb{R}^{n}\right)$.
	For any $\rho (\cdot )=(\varphi (\cdot ) ,\psi (\cdot ),\gamma  (\cdot )) \in \mathcal{N}_{\mathbb{F} }^{2}(\mathbb{T} ;\mathbb{R}^{3n})$, its norm is given by
	\begin{equation}
		\begin{aligned}
			\|\rho (\cdot)\|_{\mathcal{N}_{\mathbb{F}}^{2}(\mathbb{T}; \mathbb{R}^{3n})} := & \bigg\{ \mathbb{E}\bigg[\sum_{k=0}^{N-1} |\varphi(k)|^{2}  \\
			& \quad + \sum_{k=0}^{N-1} |\psi(k)|^{2} + \sum_{k=0}^{N-1} |\gamma(k)|^{2} \bigg] \bigg\}^{\frac{1}{2}}. \nonumber
		\end{aligned}
	\end{equation}

	$\bullet$ $\mathcal{H}(\overline{\mathbb{T}}):=\ \mathbb{R}^{n}\times L_{\mathcal{F}_{N-1} }^{2}(\Omega ;\mathbb{R}^n)\times \mathcal{N}^{2}(\mathbb{T} ;\mathbb{R}^{3n} )$. For any $(\xi  ,\eta ,\rho (\cdot ))\in \mathcal{H}(\overline{\mathbb{T}})$, its norm is given by
	{\small\begin{equation}
		\begin{aligned}
			\|(\xi  ,\eta ,\rho (\cdot ))\|_{\mathcal{H} (\overline{\mathbb{T}})}:=\left \{ \|\xi \|^{2}+\|\eta \|^{2}_{L^{2} _{\mathcal{F}_{N-1} }(\Omega ;\mathbb{R}^n )}+\|\rho(\cdot) \|^{2}_{\mathcal{N}^{2}(\mathbb{T} ;\mathbb{R}^{3n} )} \right \}^{\frac{1}{2}}.\nonumber
		\end{aligned}
	\end{equation}}

	In what follows, we shall present some basic results on the  stochastic difference equation ( S$\bigtriangleup $E) and BS$\bigtriangleup $E.

	Firstly, we study the following S$\bigtriangleup $E:

	\begin{equation}\label{eq:2.1}        		
		\left\{\begin{aligned}               	 	
				x_{k+1}  =&b\left(k, x_{k}\right)+\sigma(k, x_k) \omega_k, \\
				x_0  =&\eta , \quad k\in \mathbb{T}.
			\end{aligned}\right.
		\end{equation}

		The coefficients $(b,\sigma,\eta )$ are assumed to satisfy the following conditions:

\begin{ass}\label{ass:2.1}      
			$\eta \in \mathbb{R}^n$ and $(b,\sigma)$ are two  given random mappings
			\begin{equation}
				\begin{aligned}
					&b:\Omega\times\mathbb{T}\times  \mathbb{R}^n  \to \mathbb{R}^n ,\\
					&\sigma :\Omega\times\mathbb{T}\times  \mathbb{R}^n  \to \mathbb{R}^{n}
				\end{aligned}\nonumber
			\end{equation}
			satisfying:
			\\(i) For any $x \in \mathbb{R}^n$, $b(k,x) $ and $\sigma(k,x) $ are $\mathcal{F}_{k-1}$-measurable, $ k \in \mathbb T$;
			\\(ii)  $ b(\cdot, 0)\in L_{\mathbb{F}}^2\left(\mathbb{T}; \mathbb{R}^n\right)$  and $\sigma(\cdot, 0)\in L_{\mathbb{F}}^2\left(\mathbb{T}; \mathbb{R}^{n}\right);$
			\\(iii) The mappings $b$ and $\sigma$  are uniformly Lipschitz continuous with respect to $x$, i.e., for any $x,\bar{x}\in \mathbb{R}^n$, there exists a constant $L>0$ such that
			\begin{equation}
				|b(k,x)-b(k,\bar{x})|\le L|x-\bar{x}|	,	\quad \quad\forall k \in \mathbb T,\nonumber
			\end{equation}
			and
			\begin{equation}
				|\sigma(k,x)-\sigma (k,\bar{x})|\le L|x-\bar{x}|	, \quad \quad \forall k \in \mathbb T.\nonumber
			\end{equation}

		\end{ass}	

	\begin{lem}\label{lem:2.2}     			
	Under Assumption \ref{ass:2.1}, S$\bigtriangleup $E \eqref{eq:2.1} with coefficients $( b, \sigma, \eta )$ admits a unique solution $x(\cdot )\in L_{\mathbb{F}}^2\left(\overline{\mathbb{T}}; \mathbb{R}^n\right)$. Moreover, we have the following estimate:
			{\small\begin{equation}\label{eq:2.2}
				\begin{aligned}
					{\mathbb{E}}\left[\sum_{k=0}^N\left|x_k\right|^2\right] \leq C\left\{\left|\eta \right|^2+{\mathbb{E}}\left[\sum_{k=0}^{N-1}
					\left( |b(k, 0)|^2+\left|\sigma(k, 0)\right|^2\right)\right]\right\},
				\end{aligned}
			\end{equation}}
			where $C$ is a positive constant  only depending on $\mathbb T$ and the Lipschitz constant $L$.
			Furthermore, let $(\bar{b}, \bar{\sigma}, \bar{\eta })$ be another set of coefficients satisfying Assumption \ref{ass:2.1}, and assume that $\bar{x}(\cdot )\in  L_{\mathbb{F}}^2\left(\overline{\mathbb{T}}; \mathbb{R}^n\right)$ is a solution to S$\bigtriangleup $E\eqref{eq:2.1} corresponding  the coefficients $( \bar{b}, \bar{\sigma}, \bar{\eta })$.
			Then the following estimate holds:
			\begin{align} \label{eq:2.3}
				&\mathbb{E}\left[\sum_{k=0}^{N} |x_{k} - \bar{x}_{k}|^2\right]
				\\& \le C \mathbb{E}\bigg[|\eta - \bar{\eta}|^2  \nonumber
				+  \sum_{k=0}^{N-1} |b(k, \bar{x}_{k}) - \bar{b}(k, \bar{x}_{k})|^2  \nonumber \\
				& \qquad +  \sum_{k=0}^{N-1} |\sigma(k, \bar{x}_{k}) - \bar{\sigma}(k, \bar{x}_{k})|^2\bigg], \nonumber
			\end{align}
		where $C$ is also a positive constant which only depends  on $\mathbb T$ and the Lipschitz constant $L$.
		\end{lem}

$\mathbf{Proof.}$ For simplicity, denote
			\begin{equation}
				\left\{\begin{aligned}
					& \widehat{x}_k=x_k-\bar{x}_k,\\
					&\widehat{b}_k=b(k,x_k)-\bar{b}(k,\bar{x}_k ),\\
					&\widehat{\sigma }_k=\sigma (k,x_k)-\bar{\sigma }(k,\bar{x}_k ).
				\end{aligned}\right.\nonumber
			\end{equation}
			We will  use mathematical induction to write our proof:
			for \( N = 0 \), we have
			\[
			\mathbb{E}[|\widehat{x}_0|^2] = \mathbb{E}[|\eta  - \bar{\eta }|^2].
			\]
			This satisfies the initial condition.

			Assume that the inequality
			\begin{equation}
				\begin{aligned}
					&\mathbb{E}\left[\sum_{k=0}^{N-1} |\widehat{x}_k|^2\right]
					\\&\le  C \mathbb{E} \left[\left|\eta - \bar{\eta}\right|^2
					+ \sum_{k=0}^{N-2} \left( |b(k, \bar{x}_k) - \bar{b}(k, \bar{x}_k)|^2
					\right. \right.\\
					& \qquad	+ |\sigma(k, \bar{x}_k) - \bar{\sigma}(k, \bar{x}_k)|^2 ) \bigg]
				\end{aligned}
			\end{equation}
			holds for some constant \( C \) depending on \( \mathbb{T} \) and the Lipschitz constant \( L \).
			We need to show that

			\[
			\mathbb{E}\left[\sum_{k=0}^{N} |\widehat{x}_k|^2\right] \le C \mathbb{E} \left[\left|\eta  - \bar{\eta }\right|^2 + \sum_{k=0}^{N-1} \left( |b(k, \bar{x}_k) - \bar{b}(k, \bar{x}_k)|^2\right. \right.
			\]
			\[
			\left. \qquad \qquad \quad + |\sigma(k, \bar{x}_k) - \bar{\sigma}(k, \bar{x}_k)|^2 \right) \bigg].
			\]
			In fact, we have
			$$
			\mathbb{E}[|\widehat{x}_N|^2] = \mathbb{E}\left[|\widehat{b}_{N-1} + \widehat{\sigma}_{N-1} \omega_{N-1}|^2\right].
			$$
			Expanding this, we get
			\[
			\mathbb{E}[|\widehat{x}_N|^2] = \mathbb{E}\bigg[|\widehat{b}_{N-1}|^2 + 2 \langle \widehat{b}_{N-1}, \widehat{\sigma}_{N-1} \omega_{N-1} \rangle
			+ |\widehat{\sigma}_{N-1} \omega_{N-1}|^2 \bigg].
			\]
			Since \( \mathbb{E}[\langle \widehat{b}_{N-1}, \widehat{\sigma}_{N-1} \omega_{N-1} \rangle] = 0 \) and \( \mathbb{E}[|\widehat{\sigma}_{N-1} \omega_{N-1}|^2] \\= \mathbb{E}[|\widehat{\sigma}_{N-1}|^2] \), we obtain
			\[
			\mathbb{E}[|\widehat{x}_N|^2] \leq \mathbb{E}[|\widehat{b}_{N-1}|^2] + \mathbb{E}[|\widehat{\sigma}_{N-1}|^2].
			\]
			Using the Lipschitz condition:
			\[
\begin{split}
    &\mathbb{E}[|\widehat{b}_{N-1}|^2] \\& \leq C \mathbb{E}\bigg[|\widehat{x}_{N-1}|^2
    \quad + |b(N-1, \bar{x}_{N-1}) - \bar{b}(N-1, \bar{x}_{N-1})|^2\bigg]
\end{split}
\]
and
\[
\begin{split}
			&\mathbb{E}[|\widehat{\sigma}_{N-1}|^2]
 \\&\leq C \mathbb{E}\left[|\widehat{x}_{N-1}|^2 + |\sigma(N-1, \bar{x}_{N-1}) - \bar{\sigma}(N-1, \bar{x}_{N-1})|^2\right].
\end{split}
			\]
			Combining these results yields
			\begin{equation}
				\begin{aligned}
					&\mathbb{E}\left[\left|\widehat{x}_N\right|^2\right]
\\& \leq C \mathbb{E} \bigg[
					|\widehat{x}_{N-1}|^2 + |b(N-1, \bar{x}_{N-1}) - \bar{b}(N-1, \bar{x}_{N-1})|^2 \\
					&\quad + \left|\sigma(N-1, \bar{x}_{N-1}) - \bar{\sigma}(N-1, \bar{x}_{N-1})\right|^2
					\bigg] .
				\end{aligned}
			\end{equation}
			By the inductive hypothesis, we get
			\begin{equation}
				\begin{aligned}
					\mathbb{E}\left[\sum_{k=0}^{N} |\widehat{x}_k|^2\right] \le C \mathbb{E}
					&\bigg[|\eta  - \bar{\eta }|^2
					+ \sum_{k=0}^{N-1} \bigg( |b(k, \bar{x}_k) - \bar{b}(k, \bar{x}_k)|^2 \\
					&+ |\sigma(k, \bar{x}_k) - \bar{\sigma}(k, \bar{x}_k)|^2 \bigg) \bigg].\nonumber
				\end{aligned}
			\end{equation}
			Thus, the inductive step holds. By induction, the result is established for all \( N \).  In the above,  $C$  is a general positive constant that can be changed line by line. In particular, setting \( (\bar{b}, \bar{\sigma}, \bar{\eta }) = (0,0,0) \) confirms the result for \( \eqref{eq:2.2} \). This completes the proof.  \qed

Secondly, we consider the BS$\bigtriangleup $E as follows:
		\begin{equation}\label{eq:2.6}
			\left\{\begin{aligned}
				&y_k =f\left(k+1,y'_{k+1},z'_{k+1}\right), \\
				\\&y_N  =\xi,\quad k \in \mathbb{
					T},
			\end{aligned}\right.
		\end{equation}
		where $z_{k+1}=y_{k+1}w_k,y'_{k+1}=\mathbb{E} [y_{k+1 }|{\mathcal{F}_{k-1}}]$, $z'_{k+1}=\mathbb{E} [z_{k+1}|{\mathcal{F}_{k-1}}]$. The coefficients $(\xi, f)$ are assumed to satisfy the following assumptions:
		\begin{ass}\label{ass:2.2}
			$\xi \in L^2_{\mathcal{F}_{N-1} } (\Omega ;\mathbb{R}^n)$,  $f$ is a given  mapping:	$ f: \Omega \times \mathbb{T}  \times \mathbb{R}^n \times \mathbb{R}^{n}\to \mathbb{R}^n$
			satisfying: \\
			(i) For any $y',z' \in \mathbb{R}^n$, $f(k+1,y',z') $ is $\mathcal{F}_{k-1}$-measurable, $k \in \mathbb{
				T}$;
			\\(ii)  $ f(\cdot, 0,0)\in L_{\mathbb{F}}^2\left(\mathbb{T}; \mathbb{R}^n\right)$;
			\\(iii) the mappings $f$ is uniformly Lipschitz continuous with respect to ($y',z'$), i.e., for any $y',\bar{y}',z',\bar{z}'\in \mathbb{R}^{n}$,
			$k \in \mathbb{
				T}$, there exists a constant $L>0$ such that
			\begin{eqnarray}
				\begin{aligned}
					|f&(k+1 , y', z')-f(k+1 ,\bar{y}' , \bar{z}')|\le L\big(|y'-\bar{y}' |+|z'-\bar{z}' |\big).\nonumber
				\end{aligned}
			\end{eqnarray}
		\end{ass}
		\begin{lem}\label{lem:2.3}
			Let the coefficients $(\xi, f)$  satisfy Assumption \ref{ass:2.2}. Then  BS$\bigtriangleup $E \eqref{eq:2.6} with  the coefficients $(\xi, f)$
			admits a unique solution $y(\cdot )\in
			L_{\mathbb{F}}^2\left(\overline{\mathbb{T}}
			; \mathbb{R}^{n}\right)$. Moreover, the following estimate holds:
			\begin{equation} \label{eq:2.7}
				{\mathbb{E}}\left[\sum_{k=0}^N\left|y_k\right|^2\right] \leq C\left\{\mathbb E\left|\xi\right|^2+{\mathbb{E}}\left[\sum_{k=0}^{N-1}
				\left(|f(k+1, 0,0)|^2\right)\right]\right\},
			\end{equation}where $C$ is a positive constant  depending on $\mathbb{T}$ and the  Lipschitz constant $L$. Furthermore, suppose that $\bar{y}(\cdot )\in
			L_{\mathbb{F}}^2\left(\mathbb{T}
			; \mathbb{R}^{n}\right)$is
			a solution to BS$\bigtriangleup $E \eqref{eq:2.6} with another  coefficients $(\bar{\xi}, \bar{f})$ satisfying Assumption \ref{ass:2.2}.
			Then we have the following estimate:
			\begin{equation} \label{eq:2.8}
				\begin{aligned}
					&{\mathbb{E}}\bigg[\sum_{k=0}^N|y_k-\bar{y}_k|^2\bigg] \\
					& \leq C\bigg\{{\mathbb{E}}\bigg[\xi-\bar{\xi}|^2\\
					+&\sum_{k=0}^{N-1}
					\bigg(|f(k+1,\bar{y}'_{k+1},\bar{z}'_{k+1})-\bar{f}(k+1,\bar{y}'_{k+1},\bar{z}'_{k+1})|^2\bigg)\bigg]\bigg\}.
				\end{aligned}
			\end{equation}
			where $z_{k+1}=y_{k+1}w_k,y'_{k+1}=\mathbb{E} [y_{k+1 }|{\mathcal{F}_{k-1}}]$, $z'_{k+1}=\mathbb{E} [z_{k+1}|{\mathcal{F}_{k-1}}]$ and $C$ is a positive constant only depending on $\mathbb{T}$ and the  Lipschitz constant $L$.
		\end{lem}

$\mathbf{Proof.}$ Using Jensen's inequality and the Holder's inequality, we can derive
			\begin{equation}
				\begin{aligned}
					{\mathbb{E}}\left[\left|y_k^{\prime}\right|^2\right]  &={\mathbb{E}}\bigg[ |{\mathbb{E}}\left[y_k \mid\mathcal{F}_{k-2}\right]|^2 \bigg] \\
					&\leq {\mathbb{E}}\left[{\mathbb{E}}\left[|y_k|^2\mid \mathcal{F}_{k-2}\right]\right]
					={\mathbb{E}}\left[|y_k|^2\right]
				\end{aligned}
			\end{equation}
			and
			\begin{equation}
				\begin{aligned}
					{\mathbb{E}}\left[\left|z_k^{\prime}\right|^2\right] & ={\mathbb{E}}\bigg[|{\mathbb{E}}\left[y_k w_{k-1} \mid \mathcal{F}_{k-2}\right]|^2 \bigg]\\
					& \left.\leq {\mathbb{E}}\bigg[{\mathbb{E}}\left[|y_k|^2\mid \mathcal{F}_{k-2}\right] \cdot {\mathbb{E}}\left[w_{k-1}^2 \mid \mathcal{F}_{k-2}\right]\right] \\
					& ={\mathbb{E}}\left[{\mathbb{E}}\left[|y_k|^2 \mid \mathcal{F}_{k-2}\right]\right] \\
					& ={\mathbb{E}}\left[|y_k|^2\right].
				\end{aligned}
			\end{equation}
			Similarly, we can get an estimate for $z_k'-\bar{z}'_k$ as follows:
			\begin{equation}  \label{eq:2.11}
				\begin{aligned}
					\mathbb{E}\left[\left|z_k^{\prime}-\bar{z}'_k\right|^2\right] & =\mathbb{E}\left[|\mathbb{E}\left[\left.\left(y_k-\bar{y}_k\right) w_{k-1}\right| \mathcal{F}_{k-2} \right]|^2\right] \\
					& \leq \mathbb{E}\left[\mathbb{E}\left[|y_k-\bar{y}_k|^2 \mid \mathcal{F}_{k-2}\right] \cdot E\left[w_{k-1}^2|\mathcal{F}_{k-2}\right]\right] \\
					& =\mathbb{E}\left[\left|y_k-\bar{y}_k\right|^2\right].
				\end{aligned}
			\end{equation}
			Therefore, using the elementary inequality $|a+b|^2\le  2( |a|^2+|b|^2)$ and  the Lipschitz condition of $f$, we can derive

			\begin{equation} \label{eq:2.12}
				\begin{aligned}
					&\mathbb{E}\left[\left|y_{k-1}-\bar{y}_{k-1}\right|^2\right] \\
					&= \mathbb{E}\left[\left|f\left(k, y'_k, z'_k\right)-\bar{f}\left(k, \bar{y}'_k, \bar{z}'_k\right)\right|^2\right] \\
					&= \mathbb{E}\left[\left|f\left(k, y'_k, z'_k\right) - f\left(k, \bar{y}'_k, \bar{z}'_k\right) \right.\right. \\
					&\qquad \quad + \left. \left. f\left(k, \bar{y}'_k, \bar{z}'_k\right) - \bar{f}\left(k, \bar{y}'_k, \bar{z}'_k\right)\right|^2\right] \\
					&\leq C \cdot \mathbb{E} \left[\left|f\left(k, y'_k, z'_k\right) - f\left(k, \bar{y}'_k, \bar{z}'_k\right)\right|^2\right. \\
					&\qquad \quad + \left.\left|f\left(k, \bar{y}'_k, \bar{z}'_k\right) - \bar{f}\left(k, \bar{y}'_k, \bar{z}'_k\right)\right|^2\right] \\
					&\leq C \cdot \mathbb{E} \left\{\left[\left|y'_k - \bar{y}'_k\right| + \left|z'_k - \bar{z}'_k\right|\right]^2 \right. \\
					&\qquad \quad + \left.\left|f\left(k, \bar{y}'_k, \bar{z}'_k\right) - \bar{f}\left(k, \bar{y}'_k, \bar{z}'_k\right)\right|^2\right\} \\
					&\leq C \cdot \mathbb{E} \left\{\left|y_k - \bar{y}_k\right|^2 \right. \\
					&\qquad \quad + \left.\left|f\left(k, \bar{y}'_k, \bar{z}'_k\right) - \bar{f}\left(k, \bar{y}'_k, \bar{z}'_k\right)\right|^2\right\},
				\end{aligned}
			\end{equation}
			where the last inequality we used (\ref{eq:2.11})  and C is a general positive constant that can be changed line by line.

			Based on (\ref{eq:2.12}),  by induction we can directly obtain  \eqref{eq:2.8}. In particular, if we choose $(\bar{f},\bar{\xi})=(0,0)$, from \eqref{eq:2.8},  we can obtain  \eqref{eq:2.7}. This proof is complete. \qed

 In what follows, we present the main results of this paper.

\section{FBS$\bigtriangleup $Es with Domination-Monotonicity Conditions}\label{sec:3}
		In this section, we will be committed to studying the FBS$\bigtriangleup $E \eqref{eq:1.1}:
		\begin{equation}
			\left\{\begin{aligned}
				&x_{k+1}  =b\left(k, x_k, y'_{k+1},z'_{k+1}\right)+\sigma(k, x_k, y'_{k+1},z'_{k+1}) \omega_k ,\\
				&y_k = -f\left(k+1, x_k, y'_{k+1},z'_{k+1}\right),\\
                &z_{k+1} =y_{k+1}w_k,
				\\&y'_{k+1}=\mathbb E[y_{k+1}|\mathcal{F}_{k-1}],
				\\&z'_{k+1}=\mathbb E[z_{k+1}|\mathcal{F}_{k-1}],
				\\ &x_0  =\Lambda (y_0),
				\\&y_N  =\Phi(x_N), \quad k \in \mathbb{T}.
			\end{aligned}\right.
		\end{equation}
		Now,we  make the following assumptions on the coefficients $(\Lambda ,\Phi ,\Gamma)$ of FBS$\bigtriangleup $E \eqref{eq:1.1}.

		\begin{ass}\label{ass:3.1}
			(i)For any $x ,y,z\in \mathbb{R}^n$, $\Phi (x) \in L^{2} _{\mathcal{F}_{N-1} }(\Omega ; \mathbb {R}^n)$
			and $\Lambda (y)$ is deterministic. Furthermore, for any
			$\theta=(x,y',z')\in \mathbb{R}^{3n}$, $
			\Gamma (k,\theta )=(f(k+1,\theta) ,b(k,\theta ),\sigma (k,\theta ))
			\in \mathcal{N}^{2}({\mathbb{T}};\mathbb{R}^{3n} )$. Moreover, $(\Lambda(0),\Phi (0),\Gamma (\cdot ,0,0,0) )\in \mathcal{H}(\overline{\mathbb{T}})$;
			(ii)The mappings $\Phi$, $\Gamma$ are uniformly Lipschitz continuous, i.e., for any
			$ \theta ,\bar{\theta } \in \mathbb{R}^{3n}, k\in \mathbb T$, there exists a constant $L>0$ such that
			\begin{equation}
				\left\{\begin{aligned}
					&|\Lambda(y)-\Lambda(\bar{y})|\le L|y-\bar{y}|,\\
					&|\Phi (x)-\Phi (\bar{x} )|\le L|x-\bar{x} |,\\
					&|h(k,\theta )-h(k,\bar{\theta} )|\le L|\theta -\bar{\theta }|,\\
					&|f(k+1,\theta)-f(k+1,\bar{\theta})|\le L|\theta -\bar{\theta }|,
				\end{aligned}\right.\nonumber
			\end{equation}
			where $h=b,\sigma$.
		\end{ass}
		In addition to Assumption \ref{ass:3.1} above, we continue to introduce the following domination-monotonicity conditions on the coefficients $(\Lambda,\Phi ,\Gamma)$ for the convenience of future study.
		\begin{ass}\label{ass:3.2}
			There exist two constants $\mu\ge 0$, $v\ge0$, a matrix \(M \in \mathbb{R}^{\bar{m}\times n}\), a matrix-valued random variable $G\in L^{\infty}_{\mathcal{F}_{N-1} }(\Omega ;\mathbb{R}^{\tilde m\times n} )$, and a series of matrix-valued processes $A(\cdot ), B(\cdot ),C(\cdot )\in  L_{\mathbb{F}}^{\infty}\left(\mathbb{T}; \mathbb{R}^{m\times n}\right)$ (where \(m, \bar{m}, \tilde{m} \in \mathbb{N}\) are given)
			such that we have the following conditions:

			(i) One of the following two cases holds true. Case $1$: $\mu>0$ and $v=0$. Case $2$: $\mu=0$ and $v>0$.

			(ii) (domination condition) For  all $k\in \mathbb{T} $ and almost all $\omega \in \Omega$, and any $x,\bar{x}, y,\bar{y},z,\bar{z}\in \mathbb{R}^n$,
			\begin{equation}\label{eq:3.1}
				\left\{\begin{aligned}
					&|\Lambda(y)-\Lambda(\bar{y})|\le \frac{1}{\mu }|M\widehat{y}|,\\
					&|\Phi (x)-\Phi (\bar{x} )|\le \frac{1}{v}|G\widehat{x} |,\\
					& |f(k+1,x,y,z)-f(k+1,\bar{x},y,z)| \le \frac{1}{v}\big|A_k\widehat{x}\big|,\\
					& |h(k,x,y,z)-h(k,x,\bar{y},\bar{z})| \le \frac{1}{\mu}\bigg|B_k\widehat{y}+C_k\widehat{z}\bigg|,
				\end{aligned}\right.
			\end{equation}
			where $h=b, \sigma$, and $\widehat{x}=x-\bar{x}$,  $\widehat{y}=y-\bar{y}$, $\widehat{z}=z-\bar{z}$.

			It should be noticed that there is a little abuse of notations in the conditions above, when $\mu=0$ (resp. $v=0$), $\frac{1}{\mu} $ (resp. $\frac{1}{v} $) means $+\infty$. In other words, if $\mu=0$ or $v=0$, the corresponding domination constraints will vanish.

			(iii) (monotonicity condition) For  all $k\in \mathbb{T} $ and almost all $\omega \in \Omega$ and any $\theta =(x,y,z),\bar{\theta}=(\bar{x},\bar{y},\bar{z}) \in \mathbb{R}^{3n}$ ,
			\begin{equation}\label{eq:3.2}
				\left\{\begin{aligned}
					&\langle\Lambda(y)-\Lambda(\bar{y}), \widehat{y}\rangle \le -\mu|M \widehat{y}|^2,\\
					&\left \langle\Phi (x)-\Phi (\bar{x} ),\widehat{x}\right \rangle\ge v|G\widehat{x}|^2,\\
					& \langle\Gamma (k,\theta)-\Gamma (k,\bar{ \theta}),\widehat{\theta}\rangle
					\leq -v\big|A_k\widehat{x}\big|^2-\mu \bigg|B_k\widehat{y}+C_k\widehat{z}\bigg|^2,
				\end{aligned}\right.
			\end{equation}
			where $\langle\Gamma (k,\theta),\theta\rangle =\langle f(k+1,\theta),x \rangle +\langle b(k,\theta),y'\rangle+\langle \sigma(k,\theta),z'\rangle$
			and $\widehat{\theta } =\theta - \bar{\theta }$.
		\end{ass}
		\begin{remark}\label{rmk:3.1}
			(i) In Assumption \ref{ass:3.2}-(ii), the constant $1/\mu$ and $1/v$ can be replaced by $K/\mu$ and $K/v  (K>0)$. Therefore, for simplicity, we prefer to omit the constant $K$;\\
			(ii) There exists a symmetrical version of Assumption \ref{ass:3.2}-(iii) as follows:

			For  all $k\in \mathbb{T} $ and almost all $\omega \in \Omega$ and any $\theta,\bar{\theta} \in \mathbb{R}^{3n}$(the argument $k$ is suppressed) ,
			\begin{equation}\label{eq:3.3}
				\left\{\begin{aligned}
					&\langle\Lambda(y)-\Lambda(\bar{y}), \widehat{y}\rangle \ge \mu|M \widehat{y}|^2,\\
					&\left \langle\Phi (x)-\Phi (\bar{x} ),\widehat{x}\right \rangle\le -v|G\widehat{x}|^2,\\
					&\langle\Gamma (k,\theta)-\Gamma (k,\bar{ \theta}),\widehat{\theta}\rangle
					\ge v\big|A_k\widehat{x}\big|^2+\mu \bigg|B_k\widehat{y}+C_k\widehat{z}\bigg|^2.
				\end{aligned}\right.
			\end{equation}
			It is easy to verify the symmetry between (\ref{eq:3.2}) and (\ref{eq:3.3}), and we omit the detailed proofs. For similar proofs, we can refer to  \cite{yu2022forward}.
		\end{remark}
		For notational convenience, we would like to give some notations as follows:
		\begin{equation}\label{eq:6.4}
			\left\{\begin{aligned}
				&P(k,x_k)=A_kx_k,\\
				&P(k,\widehat{x}_k)=A_k\widehat{x}_k,\\
				&Q(k,y'_{k+1},z'_{k+1})=B_ky'_{k+1}+C_kz'_{k+1},\\
				&Q(k,\widehat{y}'_{k+1},\widehat{z}'_{k+1})=B_k\widehat{y}'_{k+1}+C_k\widehat{z}'_{k+1} ,\quad k \in\mathbb{T}.
			\end{aligned}\right.
		\end{equation}

		Now, we present the main results of this section.
		\begin{thm}\label{thm:3.1}
			Let $(\Lambda,\Phi ,\Gamma)$ be a set of coefficients satisfying Assumption \ref{ass:3.1} and Assumption \ref{ass:3.2}. Then FBS$\bigtriangleup $E \eqref{eq:1.1} admits a unique solution $(x(\cdot),y(\cdot))\in N^{2}(\overline{\mathbb{T}};\mathbb{R}^{2n} )$. Moreover, we have the following estimate:
			\begin{equation}\label{eq:3.4}
				\mathbb E\bigg[\displaystyle \sum_{k=0}^N|{x}_k|^2+\displaystyle\sum_{k=0}^N|{y}_k|^2\bigg ]\le K\mathbb{E}[\mathrm{I}],
			\end{equation}
			where
			\begin{equation}
				\begin{aligned}
					\mathrm{I}=&|\Phi (0)|^2+\sum_{k=0}^{N-1}|b(k,0,0,0)|^2 +\sum_{k=0}^{N-1}|\sigma (k,0,0,0)|^2\\
					&+\sum_{k=0}^{N-1}|f(k+1,0,0,0)|^2+|\Lambda (0)|^2 ,
				\end{aligned}
			\end{equation}
			and $K$ is a positive constant depending only on $\mathbb T$, the Lipschitz constants, $\mu$, $v$ and the bounds of all $G$, $A(\cdot)$,  $B(\cdot)$, $C(\cdot), M$. Furthermore, let $(\bar{\Lambda},\bar{ \Phi },\bar{ \Gamma})$ be another set of coefficients satisfying Assumption \ref{ass:3.1} and Assumption \ref{ass:3.2}, and $ (\bar{x}(\cdot),\bar{y}(\cdot))\in N^{2}(\overline{\mathbb{T}};\mathbb{R}^{2n} )$ be a solution to FBS$\bigtriangleup $E \eqref{eq:1.1} with the coefficients $(\bar{\Lambda},\bar{ \Phi },\bar{ \Gamma})$.  Then the following estimate holds:
			\begin{equation}\label{eq:3.6}
				\mathbb E\bigg[\displaystyle \sum_{k=0}^N|\widehat{x}_k|^2+\displaystyle\sum_{k=0}^N|\widehat{y}_k|^2\bigg ]\le K\mathbb{E}[\widehat{\mathrm{I}}],
\end{equation}
where
			\begin{equation}
				\begin{aligned}
					\widehat{\mathrm{I}}=&|\Phi (\bar{x}_N)-\bar{\Phi} (\bar{x}_N)|^2+\sum_{k=0}^{N-1}|b(k,\bar{\theta}_k)-\bar{b}(k,\bar{\theta}_k)|^2\\
					&+\sum_{k=0}^{N-1}|\sigma(k,\bar{\theta}_k)-\bar{\sigma}(k,\bar{\theta}_k)|^2
					\\&+\sum_{k=0}^{N-1}|f(k+1,\bar\theta_k)-\bar{f}(k+1,\bar{\theta}_k)|^2+|\Lambda (\bar{y}_0)-\bar \Lambda(\bar{y}_0)|^2
				\end{aligned}
			\end{equation}
			and \begin{equation*}
				\begin{aligned}
					\theta_k &= \left(x_k, y_{k+1}'p, z_{k+1}'\right) \\
					&= \left(x_k, \mathbb{E}\left[y_{k+1} \mid \mathcal{F}_{k-1}\right], \mathbb{E}\left[y_{k+1} w_k \mid \mathcal{F}_{k-1}\right]\right), \\
					\bar{\theta}_k &= \left(\bar{x}_k, \bar{y}_{k+1}', \bar{z}_{k+1}'\right).
				\end{aligned}
			\end{equation*}
In addition, $K$ is the same constant as in \eqref{eq:3.4}.
		\end{thm}
		Next, we are devoted to proving Theorem \ref{thm:3.1}. Due to the symmetry of monotonicity conditions \eqref{eq:3.2} and \eqref{eq:3.3}, we only give the detailed proofs under the monotonicity condition \eqref{eq:3.2}.

		For any $\big(\xi ,\eta ,\rho (\cdot )\big)\in \mathcal{H} (\overline{\mathbb{T}})$ with $\rho (\cdot )=(\varphi (\cdot ) ,\psi (\cdot ),\gamma  (\cdot ))$  we continue to introduce a family of FBS$\bigtriangleup $Es parameterized by $\alpha\in[0,1]$ as follows:
		\begin{equation}\label{eq:3.8}
			\left\{\begin{aligned}
				x^{\alpha}_{k+1} =&\big[b^{\alpha}(k,\theta_{k}^{\alpha}) +\psi_k\big]+\big[\sigma^{\alpha}(k,\theta_{k}^{\alpha}) +\gamma_k\big] \omega_k,    \\
				y_k^{\alpha}=& -\big[f^{\alpha}(k+1,\theta_k^{\alpha})+\varphi_k\big], \quad k \in\mathbb{T},\\
				x_0^{\alpha} =&\Lambda^\alpha(y_0^{\alpha})+\xi,\\
				y_N^{\alpha} =&\Phi^{\alpha} (x_N^{\alpha})+\eta,
			\end{aligned}\right.
		\end{equation}
		where \begin{equation*}
			\begin{aligned}
				\theta^\alpha_k &= \left({x^\alpha_k}, {y'^\alpha_{k+1}}, {z'^\alpha_{k+1}}\right) \\
				&= \left(x^{\alpha \top}_k, \mathbb{E}\left[y^\alpha_{k+1} \mid \mathcal{F}_{k-1}\right]^{\top}, \mathbb{E}\left[y^\alpha_{k+1} \omega_k \mid \mathcal{F}_{k-1}\right]^{\top}\right)^{\top}.
			\end{aligned}
		\end{equation*}
		Here, we have used the following notation: for any $(k,\omega,\theta)\in\mathbb{T}\times\Omega\times N^{2}(\overline{\mathbb{T}};\mathbb{R}^{3n} )$,\\
		\begin{equation}\label{eq:3.9}
			\left\{\begin{aligned}
				&\Phi ^\alpha (x)=\alpha \Phi (x)+(1-\alpha)vG^{\top}Gx,\\
				&\Lambda^{\alpha}(y)=\alpha\Lambda(y)-(1-\alpha)\mu M^\top My,\\
				& f^\alpha (k+1,\theta)=\alpha f(k+1,\theta )-(1-\alpha)v\bigg[A_k^{\top}P(k,x_k)\bigg],\\
				& b^\alpha (k,\theta)=\alpha b(k,\theta)-(1-\alpha)\mu\bigg[ B_k^\top Q(k,y'_{k+1},z'_{k+1})\bigg],\\
				& \sigma ^\alpha (k,\theta)=\alpha \sigma (k,\theta )-(1-\alpha)\mu \bigg[C_k^\top Q(k,y'_{k+1},z'_{k+1}) \bigg],\\
			\end{aligned}\right.
		\end{equation}
		where $P(k,x_k)$ ,$Q(k,y'_{k+1},z'_{k+1})$ are defined by \eqref{eq:6.4}.
		Similarly, we continue to denote $\Gamma ^\alpha (k,\theta ):=(f^\alpha(k+1,\theta ) ,b^\alpha(k,\theta ) ,\sigma ^\alpha
		(k,\theta ) )$.
		Without loss of generality, we assume that the Lipschitz constants of the coefficients $(\Lambda,\Phi,\Gamma)$ are larger than
		\begin{equation}
			\begin{aligned}
				\max \{\mu, v\}\bigg(&|G|+\|A(\cdot)\|_{L_{\mathbb{F}}^{\infty}(\mathbb{T}; \mathbb{R}^{m\times n})}+\|B(\cdot)\|_{L_{\mathbb{F}}^{\infty}(\mathbb{T}; \mathbb{R}^{m\times n})}\\
				&+\|C(\cdot)\|_{L_{\mathbb{F}}^{\infty}(\mathbb{T}; \mathbb{R}^{m\times n})}\bigg)^2\nonumber
			\end{aligned}
		\end{equation}
		and the constant $\mu$ and $v$ in Assumption \ref{ass:3.2}-(i) satisfy the following condition:
		\begin{equation}
			\begin{aligned}
				\mbox{$(\frac{1}{\mu })^2,(\frac{1}{v})^2$} \ge \text{max}\Big\{&|G|,\|A(\cdot)\|_{L_{\mathbb{F}}^{\infty}(\mathbb{T}; \mathbb{R}^{m\times n})},\|B(\cdot)\|_{L_{\mathbb{F}}^{\infty}(\mathbb{T}; \mathbb{R}^{m\times n})},\\
				&\|C(\cdot)\|_{L_{\mathbb{F}}^{\infty}(\mathbb{T}; \mathbb{R}^{m\times n})}\Big\}.\nonumber
			\end{aligned}
		\end{equation}
		Then we can easily verify that for any $\alpha\in[0,1]$, the new coefficients $(\Lambda^\alpha,\Phi^\alpha,\Gamma^\alpha)$ also satisfy Assumption \ref{ass:3.1} and Assumption \ref{ass:3.2}  with the same Lipschitz constants, $\mu$, $v$, $G$, $A(\cdot)$,  $B(\cdot)$, $C(\cdot)$ as the original coefficients $(\Lambda,\Phi ,\Gamma)$.

		Obviously, when $\alpha=0$, FBS$\bigtriangleup $E \eqref{eq:3.8} can be rewritten in the following form:
		\begin{equation}\label{eq:3.10}
			\left\{\begin{aligned}
				x_{k+1}
				^{0}=&\bigg \{-\mu \Big[B_k^\top  Q(k,{y'}_{k+1},{z'}_{k+1})\Big]+\psi_k\bigg \}\\
				&+\bigg \{-\mu\Big[ C(k)^\top Q(k,{y'}_{k+1},{z'}_{k+1})\Big]+\gamma_k \bigg\}\omega_k,\\
				y_{k}^0=&-\bigg\{-v\Big[A_k^{\top}P(k,x_k)\Big]+\varphi_k\bigg\},\quad k\in \mathbb{T},\\
				x_0^{0} =& -\mu M^\top M y_0+\xi,\\
				y_N^0=&\: vG^{\top}Gx_N+\eta.
			\end{aligned}\right.
		\end{equation}
		At this point, it is easy to see that  FBS$\bigtriangleup $E\eqref{eq:3.10} is in a decoupled form.  Furthermore, when Assumption \ref{ass:3.2}-(i)-Case 1 holds (i.e., $\mu>0$ and $v=0$), we can firstly solve $y^0(\cdot)$ from the backward equation (\ref{eq:3.10}), then substitute $y^0(\cdot)$ into the forward equation (\ref{eq:3.10}) and solve $x^0(\cdot)$. Similarly, when Assumption \ref{ass:3.2}-(i)-Case 2 holds (i.e., $\mu=0$ and $v>0$), we can firstly solve the forward equation and then the backward equation. In short, when $\alpha=0$, under Assumptions \ref{ass:3.1} and \ref{ass:3.2}, FBS$\bigtriangleup $E \eqref{eq:3.10} admits a unique solution $(x^0(\cdot),y^0(\cdot))\in N^2(\overline{\mathbb{T}};\mathbb{R}^{2n} )$.

		It is clear that when $\alpha=1$ and $(\xi ,\eta ,\rho (\cdot ))$ vanish, FBS$\bigtriangleup $E \eqref{eq:3.8} and FBS$\bigtriangleup $E \eqref{eq:1.1} are identical.  Next, we will illustrate that if for some $\alpha_{0}\in[0,1)$, FBS$\bigtriangleup $E\eqref{eq:3.8} is uniquely solvable for any $\big(\xi ,\eta ,\rho (\cdot )\big)\in \mathcal{H} (\overline{\mathbb{T}})$, then there exists a fixed step length $\delta_0>0$ such that the same conclusion still holds for any $\alpha\in[\alpha_{0},\alpha_{0}+\delta_0]$. As long as this has been proved to be true, we can gradually increase the parameter $\alpha$ until $\alpha=1$. This method is called the method of continuation which was initially introduced by \cite{hu1995solution}.

		For this goal, we shall first establish an a priori estimate for the solution of FBS$\bigtriangleup $E \eqref{eq:3.8}, which plays an important role in the subsequent proofs.
		\begin{lem}\label{lem:3.2}
			Let the given coefficients $(\Lambda,\Phi ,\Gamma)$ satisfy Assumption \ref{ass:3.1} and Assumption \ref{ass:3.2}. Let $\alpha\in[0,1]$, $\big(\xi ,\eta ,\rho (\cdot )\big)$, $\big(\bar{
				\xi} ,\bar{\eta} ,\bar{\rho} (\cdot )\big)\in \mathcal{H} (\overline{\mathbb{T}})$. Assume that $(x(\cdot),y(\cdot))\in N^2(\overline{\mathbb{T}};\mathbb{R}^{2n} )$ is the solution of FBS$\bigtriangleup $E \eqref{eq:3.8} with the coefficients $(\Lambda^\alpha+\xi,\Phi^\alpha+\eta ,\Gamma^\alpha+\rho)$ and $(\bar{x}(\cdot),\bar{y}(\cdot))\in N^2(\overline{\mathbb{T}};\mathbb{R}^{2n} )$ is also the solution to FBS$\bigtriangleup $E \eqref{eq:3.8} with the coefficients $(\Lambda^\alpha+\bar{\xi},\Phi^\alpha+\bar{\eta} ,\Gamma^\alpha+\bar{\rho})$. Then the following estimate holds:
			\begin{equation}\label{eq:3.11}
				\mathbb E\bigg[\displaystyle \sum_{k=0}^N |\widehat{x}_k|^2+\displaystyle \sum_{k=0}^N |\widehat{y}_k|^2\bigg ]\le K\mathbb{E}[\widehat{\mathrm{J}}],
			\end{equation}
			where
			\begin{equation}\label{eq:3.12}
				\begin{aligned}
					\widehat{\mathrm{J} }=&|\widehat{\eta } |^2+\sum_{k=0}^{N-1}|\widehat{\varphi }_k |^2+\sum_{k=0}^{N-1}|\widehat{\psi }_k |^2
					+\sum_{k=0}^{N-1}|\widehat{\gamma}_k |^2+ |\widehat{\xi}|^2
				\end{aligned}
			\end{equation}
			and $\widehat{\xi }=\xi -\bar{\xi }$, $\widehat{\varphi }=\varphi -\bar{\varphi }$, etc.
			Here  $K$ is a positive constant depending on $\mathbb T$.
		\end{lem}
		$\mathbf{{Proof.}}$
			In the following proofs, the argument $k$ is suppressed for simplicity. Besides, it should be noted that the positive constant $K$ could be changed line by line.
			By the estimate \eqref{eq:2.3} in Lemma \ref{lem:2.2},  an estimate of $\mathbb{E}\bigg[\displaystyle\sum_{k=0}^{N} |\widehat{x}_k|^2\bigg ]$ can be obtained:
			\begin{equation}\label{eq:3.13}
				\begin{aligned}
					& \mathbb{E}\bigg[\sum_{k=0}^{N} |\widehat{x}_k|^2\bigg] \\
					& \le K \mathbb{E}\Bigg\{
					\big|\alpha \big(\Lambda(y_0) - \Lambda(\bar{y}_0)\big) - (1-\alpha)\mu M^\top M\widehat{y}_0 + \widehat{\xi}\big|^2 \\
					& \qquad + \sum_{k=0}^{N-1} \bigg| \alpha \bigg(b(k, \bar{x}_k, y_{k+1}', z'_{k+1}) - b(k, \bar{x}_k, \bar{y}_{k+1}', \bar{z}_{k+1}')\bigg) \\
					& \qquad \quad - (1-\alpha)\mu \Big[B_k^\top Q(k, \widehat{y}_{k+1}', \widehat{z}_{k+1})\Big] + \widehat{\psi}_k \bigg|^2 \\
					& \qquad + \sum_{k=0}^N \Bigg| \alpha \bigg(\sigma(k, \bar{x}_k, y_{k+1}', z_{k+1}) - \sigma(k, \bar{x}_k, \bar{y}_{k+1}', \bar{z}_{k+1}')\bigg) \\
					& \qquad \quad - (1-\alpha)\mu\Big[C_k^\top Q(k, \widehat{y}_{k+1}', \widehat{z}_{k+1}')\Big] + \widehat{\gamma}_k \Bigg|^2 \Bigg\}.
				\end{aligned}
			\end{equation}
			Similarly, by applying the estimate \eqref{eq:2.6} in Lemma \ref{lem:2.3},  we can obtain an estimate of $\mathbb E\bigg[\displaystyle \sum_{k=0}^N |\widehat{y}_k|^2\bigg ]$:
			\begin{equation}\label{eq:3.14}
				\begin{aligned}
					&\mathbb E\bigg[\displaystyle \sum_{k=0}^N |\widehat{y}_k|^2\bigg ]\\
					\le & K\mathbb E\Bigg\{\Big |\alpha \big (\Phi (x_N)-\Phi (\bar{x }_N )\big )+(1-\alpha )vG^\top G\hat{x}_N+\widehat{\eta }  \Big |^2\\
					& \qquad+\sum_{k=0}^{N-1}\bigg |\alpha \bigg(f(k+1, x_k,\bar{y}_{k+1}',\bar{z}_{k+1}')
					\\& \qquad -f(k+1, \bar{x}_k,\bar{y}_{k+1}',\bar{z}_{k+1}')\bigg)\\
					&  \qquad -(1-\alpha )v\Big[A_k^\top P(k,\widehat{x}_k)\Big]+\widehat{\varphi }_k\bigg |^2\Bigg\}.
				\end{aligned}
			\end{equation}

			We consider the equation
			\begin{equation}  \label{eq:3.171}
				\begin{aligned}
					&\mathbb{E}\left[\left\langle\widehat{x}_N, \widehat{y}_N\right\rangle-\left\langle\widehat{x}_0, \widehat{y}_0\right\rangle\right]
					\\&=\displaystyle \sum_{k=0}^{N-1}\mathbb E\bigg[\left \langle \widehat{x}_{k+1},\widehat{y}_{k+1}
					\right \rangle-\left \langle \widehat{x}_k,\widehat{y}_k \right \rangle\bigg],
				\end{aligned}
			\end{equation}
			where the left side is equal

			\begin{equation} \label{eq:3.181}
				\begin{aligned}
					&\mathbb{E}\left[\left\langle\widehat{x}_N, \widehat{y}_N\right\rangle - \left\langle\widehat{x}_0, \widehat{y}_0\right\rangle\right] \\
					= & \mathbb{E}\left[\left\langle\widehat{x}_N, \widehat{\Phi}^\alpha\left(x_N^\alpha\right) + \widehat{\eta}\right\rangle - \left\langle\widehat{\Lambda}^\alpha\left(y_0^\alpha\right) + \widehat{\xi},  \widehat{y}_0\right\rangle\right] \\
					= & \mathbb{E}\left[\left\langle\widehat{x}_N, \alpha\left(\Phi\left(x_N\right) - \Phi\left(\bar{x}_N\right)\right) + (1-\alpha) v G^{\top} G \widehat{x}_N + \widehat{\eta}\right\rangle \right. \\
					&  \left. - \left\langle\alpha\left(\Lambda\left(y_0\right) - \Lambda\left(\bar{y}_0\right)\right) - (1-\alpha) \mu M^{\top} M \widehat{y}_0 + \widehat{\xi}, \widehat{y}_0\right\rangle\right] \\
					= & \mathbb{E}\left[(1-\alpha) v \left|G \widehat{x}_N\right|^2 + (1-\alpha) \mu \left|M \widehat{y}_0\right|^2 \right. \\
					& \quad + \left\langle\widehat{x}_N, \alpha\left(\Phi\left(x_N\right) - \Phi\left(\bar{x}_N\right)\right)\right\rangle \\
					& \quad \left. - \left\langle\alpha\left(\Lambda\left(y_0\right) - \Lambda\left(\bar{y}_0\right)\right), \widehat{y}_0\right\rangle + \left\langle\widehat{x}_N, \widehat{\eta}\right\rangle - \left\langle\widehat{\xi}, \widehat{y}_0\right\rangle\right]
				\end{aligned}
			\end{equation}
			and the right side is equal
			\begin{equation}\label{eq:3.191}
				\begin{aligned}
					& \sum_{k=0}^{N-1} \mathbb{E}\bigg[\left\langle \widehat{x}_{k+1}, \widehat{y}_{k+1}\right\rangle - \left\langle \widehat{x}_k, \widehat{y}_k\right\rangle\bigg] \\
					= & \sum_{k=0}^{N-1} \mathbb{E}\bigg[\mathbb{E}\left[\left\langle\widehat{x}_{k+1}, \widehat{y}_{k+1}\right\rangle \mid \mathcal{F}_{k-1}\right] - \left\langle\widehat{x}_k, \widehat{y}_k\right\rangle\bigg]\\
						= & \sum_{k=0}^{N-1} \mathbb{E}\bigg[\mathbb{E}\left[\left\langle\widehat{b}_k^\alpha + \widehat{\psi}_k + \left[\widehat{\sigma}_k^\alpha + \widehat{\gamma}_k\right] \omega_k, \widehat{y}_{k+1}\right\rangle \mid \mathcal{F}_{k-1}\right]\\
					& \qquad  - \left\langle\widehat{x}_k, \widehat{y}_k\right\rangle\bigg]\\
					= & \sum_{k=0}^{N-1} \mathbb{E}\bigg[\left\langle\widehat{b}_k^\alpha, \widehat{y}_{k+1}^{\prime}\right\rangle + \left\langle\widehat{\psi}_k, \widehat{y}_{k+1}^{\prime}\right\rangle + \left\langle\widehat{\sigma}_k^\alpha, \widehat{z}_{k+1}^{\prime}\right\rangle
\\&\qquad + \left\langle\widehat{\gamma}_k, \widehat{z}_{k+1}^{\prime}\right\rangle
		+ \left\langle\widehat{f}_{k+1}^\alpha, \widehat{x}_k\right\rangle + \left\langle\widehat{x}_k, \widehat{\varphi}_k\right\rangle\bigg]
			\\	= & \sum_{k=0}^{N-1} \mathbb{E}\bigg[\left\langle\alpha \widehat{b}_k - (1-\alpha) \mu\left[B_k^{\top} Q(k, \widehat{y}_{k+1}^{\prime}, \widehat{z}_{k+1}^{\prime})\right], \widehat{y}_{k+1}^{\prime}\right\rangle\\
			& \qquad + \left\langle\alpha \widehat{\sigma}_k - (1-\alpha) \mu\left[C_k^{\top} Q(k, \widehat{y}_{k+1}^{\prime}, \widehat{z}_{k+1}^{\prime})\right], \widehat{z}_{k+1}^{\prime}\right\rangle \\
			& \qquad + \left\langle\alpha \widehat{f}_{k+1} - (1-\alpha) v\left[A_k^{\top} P(k, \widehat{x}_k)\right], \widehat{x}_k\right\rangle \\
			& \qquad + \left\langle\widehat{\psi}_k, \widehat{y}_{k+1}^{\prime}\right\rangle + \left\langle\widehat{\gamma}_k, \widehat{z}_{k+1}^{\prime}\right\rangle + \left\langle\widehat{x}_k, \widehat{\varphi}_k\right\rangle\bigg] \\
			= & \sum_{k=0}^{N-1} \mathbb{E}\bigg[\left\langle\alpha \widehat{f}_{k+1}, \widehat{x}_k\right\rangle + \left\langle\alpha \widehat{b}_k, \widehat{y}_{k+1}^{\prime}\right\rangle + \left\langle\alpha \widehat{\sigma}_k, \widehat{z}_{k+1}^{\prime}\right\rangle \\
			& \qquad \quad- \left\langle(1-\alpha) \mu Q(k, \widehat{y}_{k+1}^{\prime}, \widehat{z}_{k+1}^{\prime}), B_k \widehat{y}_{k+1}^{\prime} + C_k \widehat{z}_{k+1}^{\prime}\right\rangle \\
			& \qquad \quad - \left\langle(1-\alpha) v P(k, \widehat{x}_k), A_k \widehat{x}_k\right\rangle \\
			& \qquad \quad+ \left\langle\widehat{\psi}_k, \widehat{y}_{k+1}^{\prime}\right\rangle + \left\langle\widehat{\gamma}_k, \widehat{z}_{k+1}^{\prime}\right\rangle + \left\langle\widehat{x}_k, \widehat{\varphi}_k\right\rangle\bigg]
\\
		= & \sum_{k=0}^{N-1} \mathbb{E}\bigg[\alpha\langle\Gamma(\theta_k) - \Gamma(\bar{\theta}_k), \widehat{\theta}_k\rangle \\
			& - (1-\alpha) \mu \left|Q(k, \widehat{y}_{k+1}^{\prime}, \widehat{z}_{k+1}^{\prime})\right|^2 - (1-\alpha) v \left|P(k, \widehat{x}_k)\right|^2 \\
			& + \left\langle\widehat{\psi}_k, \widehat{y}_{k+1}^{\prime}\right\rangle + \left\langle\widehat{\gamma}_k, \widehat{z}_{k+1}^{\prime}\right\rangle + \left\langle\widehat{x}_k, \widehat{\varphi}_k\right\rangle\bigg],
				\end{aligned}
			\end{equation}
			where $\widehat{b}_k=b(k,\theta_k) - b(k,\bar{\theta}_k)$, and similarly for $\widehat{\sigma}, \widehat{f}, \widehat{\Gamma}$, etc.
			Putting \eqref{eq:3.181} and \eqref{eq:3.191} into \eqref{eq:3.171}, we get that
			\begin{equation}\label{eq:3.15}
				\begin{aligned}
					&\mathbb{E}\Bigg \{(1-\alpha)v|G\widehat{x}_N |^2
					+\alpha \left \langle \Phi (x_N)-\Phi (\bar{x}_N),\widehat{x}_{N}  \right \rangle\\
					&\quad+(1-\alpha)\mu|M\widehat{y}_0 |^2
					-\alpha \left \langle \Lambda (y_0)-\Lambda (\bar{y}_0),\widehat{y}_{0}  \right \rangle\\
					&\quad+
					(1-\alpha )\mu \sum_{k=0}^{N-1}\big |Q(k,\widehat y_{k+1}',\widehat z'_{k+1})\big |^2
					\\&\quad + (1-\alpha )v\sum_{k=0}^{N-1}\big|P(k,\widehat{x}_k)\big|^2\\
					&\quad -\alpha \sum_{k=0}^{N-1}\langle \Gamma (k,\theta )
					-\Gamma(k,\bar{\theta} ) ,\widehat{\theta}  \rangle \Bigg \}\\
					&=\mathbb{E}\Bigg \{
					-\left \langle \widehat{\eta},\widehat{x}_N  \right \rangle +\left \langle \widehat{\xi },\widehat{y}_0  \right \rangle\\
					&\qquad \quad +\sum_{k=0}^{N-1}\bigg[ \left\langle\widehat{\psi}_k, \widehat{y}_{k+1}^{\prime}\right\rangle
					+\left\langle \widehat{\gamma}_k, \widehat{z}_{k+1}^{\prime}\right\rangle
					+\left\langle\widehat{x}_k, \widehat{\varphi}_k\right\rangle\bigg]
					\Bigg \}.
				\end{aligned}
			\end{equation}
			Therefore, combining with the monotonicity conditions in Assumption \ref{ass:3.2}-(iii), \eqref{eq:3.15} is reduced to
			\begin{equation}\label{eq:3.16}
				\begin{aligned}
					&\mathbb{E}\Bigg \{ v|G\widehat{x}_N |^2+\mu |M\widehat{y}_0|^2+v\sum_{k=0}^{N-1}\big|P(k,\widehat{x}_k)\big|^2\\
					&\qquad +\mu \sum_{k=0}^{N-1}\big |Q(k,\widehat y_{k+1}',\widehat z_{k+1}')\big |^2\Bigg \}\\
					&\le \mathbb{E}\Bigg \{-\left \langle \widehat{\eta},\widehat{x}_N  \right \rangle
					+\left \langle \widehat{\xi },\widehat{y}_0  \right \rangle\\
					&\qquad  +\sum_{k=0}^{N-1}\Big [\left\langle\widehat{\psi}_k, \widehat{y}_{k+1}^{\prime}\right\rangle
					+\left\langle \widehat{\gamma}_k, \widehat{z}_{k+1}^{\prime}\right\rangle
					+\left\langle\widehat{x}_k, \widehat{\varphi}_k\right\rangle \Big ]\Bigg \}.
				\end{aligned}
 			\end{equation}
			The following proofs will be divided into two cases according to Assumption \ref{ass:3.2}-(i).\\
			\textbf{Case 1}: $\mu>0$ and $v=0$. By applying the domination conditions \eqref{eq:3.1} in Assumption \ref{ass:3.2}-(ii) to the estimate \eqref{eq:3.13}, we get
			\begin{equation}\label{eq:6.20}
				\begin{aligned}
					\mathbb E\bigg[\displaystyle \sum_{k=0}^N |\widehat{x}_k|^2\bigg ]\le  &K\mathbb E\Bigg \{|\widehat{\xi } |^2+|M\widehat{y}_0|^2+\sum_{k=0}^{N-1}|\widehat{\psi }_k |^2+\sum_{k=0}^{N-1}|\widehat{\gamma} _k|^2\\&+\sum_{k=0}^{N-1}\big |Q(k,\widehat y_{k+1}',\widehat z_{k+1}')\big |^2\Bigg \}.
				\end{aligned}
			\end{equation}
			Applying the Lipschitz condition to the estimate \eqref{eq:3.14} leads to
			\begin{equation}\label{eq:6.22}
				\begin{aligned}
					\mathbb E\bigg[\displaystyle \sum_{k=0}^N |\widehat{y}_k|^2\bigg ]\le K\mathbb{E}\bigg \{|\widehat{\eta}|^2+\sum_{k=0}^{N-1}|\widehat{\varphi}_k|^2
					+\displaystyle \sum_{k=0}^{N} |\widehat{x}_k |^2\bigg\}.
				\end{aligned}
			\end{equation}
			Hence, combining \eqref{eq:6.20} and \eqref{eq:6.22}, it  yields
			\begin{equation}\label{eq:3.19}
				\begin{aligned}
					&\mathbb E\bigg[\displaystyle \sum_{k=0}^N |\widehat{x}_k|^2+\displaystyle \sum_{k=0}^{N} |\widehat{y}_k|^2\bigg ]\\
					&\le K\mathbb{E}\Bigg \{\widehat{\mathrm{J}}+\sum_{k=0}^N\big |Q(k,\widehat y_{k+1}',\widehat z_{k+1}')\big |^2+|M\widehat{y_0}|^2\Bigg \},
				\end{aligned}
			\end{equation}
			where $\widehat{\mathrm{J}}$ is defined by \eqref{eq:3.12}. Finally, we continue to combine \eqref{eq:3.16} and \eqref{eq:3.19} and with the inequality $ab\le\frac{1}{4\varepsilon}a^2+\varepsilon b^2$ to get
			{\small\begin{equation}\label{eq:3.20}
					\begin{aligned}
						&\quad\mathbb E\bigg[\displaystyle \sum_{k=0}^N |\widehat{x}_k|^2+\displaystyle \sum_{k=0}^N |\widehat{y}_k|^2\bigg ]\\
						&\le K\mathbb{E}\Bigg \{\widehat{\mathrm{J}}-\left \langle \widehat{\eta},\widehat{x}_N  \right \rangle
						+\left \langle \widehat{\xi },\widehat{y}_0  \right \rangle\\
						&\qquad +\sum_{k=0}^{N-1}\Big [ \left\langle\widehat{\psi}_k, \widehat{y}_{k+1}^{\prime}\right\rangle
						+\left\langle \widehat{\gamma}_k, \widehat{z}_{k+1}^{\prime}\right\rangle
						+\left\langle\widehat{x}_k, \widehat{\varphi}_k\right\rangle\big]\Bigg \}\\
						&\le K\mathbb{E}\Bigg \{\widehat{\mathrm{J}}+2\varepsilon\bigg [ \displaystyle \sum_{k=0}^{N} |\widehat{x}_k|^2+\displaystyle \sum_{k=0}^{N} |\widehat{y}_k|^2\bigg ]\Bigg\}.
					\end{aligned}
			\end{equation}}
			By taking $\varepsilon$ small enough such that $2K\varepsilon<<1$, we can easily obtain the desired estimate \eqref{eq:3.11}, and the proof in this case is finished.\\
			\textbf{Case 2}: $\mu=0$ and $v>0$. Differently, we apply the Lipschitz conditions  to the estimate \eqref{eq:3.13} to get
			\begin{equation}\label{eq:6.26}
				\begin{aligned}
					&\mathbb E\bigg[ \sum_{k=0}^N |\widehat{x}_k|^2\bigg ]\\
					\le& K\mathbb{E}\Bigg \{\displaystyle|\widehat{\xi}|^2
					+ \sum_{k=0}^{N-1} |\widehat{y}_k|^2+\sum_{k=0}^{N-1}|\widehat{\psi }_k |^2+\sum_{k=0}^{N-1}|\widehat{\gamma}_k|^2\Bigg \}.
				\end{aligned}
			\end{equation}
			By the domination conditions \eqref{eq:3.1} in Assumption \ref{ass:3.2}-(ii), we deduce from \eqref{eq:3.14} and obtain
			\begin{equation}\label{eq:3.22}
				\begin{aligned}
					\quad\mathbb E\bigg[\displaystyle \sum_{k=0}^N |\widehat{y}_k|^2\bigg ]
					&\le K\mathbb{E}\bigg \{ |\widehat{\eta}|^2
					+|G\widehat{x}_N|^2\\&+\sum_{k=0}^{N-1} |P(k,\widehat{x}_k)|^2+\sum_{k=0}^{N-1}|\widehat{\varphi}_k|^2
					\Bigg \}.
				\end{aligned}
			\end{equation}
			Thus, combining \eqref{eq:6.26} and \eqref{eq:3.22}, we can derive
			\begin{equation}\label{eq:3.23}
				\begin{aligned}
					\mathbb E\bigg[\displaystyle \sum_{k=0}^{N} |\widehat{x}_k|^2+\displaystyle \sum_{k=0}^N |\widehat{y}_k|^2\bigg ]
					\le K\mathbb{E}\Bigg \{\widehat{\mathrm{J} }+|G\widehat{x}_N|^2+\sum_{k=0}^{N-1}|P(k,\widehat{x})|^2\Bigg \},
				\end{aligned}
			\end{equation}
			where $\widehat{\mathrm{J}}$ is defined by \eqref{eq:3.12}. At last, \eqref{eq:3.16} and \eqref{eq:3.23} work together to turn out
			\begin{equation}
				\begin{aligned}
					&\quad\mathbb E\bigg[\displaystyle \sum_{k=0}^N |\widehat{x}_k|^2+\displaystyle \sum_{k=0}^N |\widehat{y}_k|^2\bigg ]\\
					&\le K\mathbb{E}\Bigg \{\widehat{\mathrm{J} }-\left \langle \widehat{\eta},\widehat{x}_N  \right \rangle +\left \langle \widehat{\xi },\widehat{y}_0  \right \rangle
					\\&+\sum_{k=0}^{N-1}\Big [ \left\langle\widehat{\psi}_k, \widehat{y}_{k+1}^{\prime}\right\rangle
					+\left\langle \widehat{\gamma}_k, \widehat{z}_{k+1}^{\prime}\right\rangle
					+\left\langle\widehat{x}_k, \widehat{\varphi}_k\right\rangle\big]\Bigg \}.
				\end{aligned}
			\end{equation}
			The remaining proof is the same as \eqref{eq:3.20} in Case 1, and then we finish the proof in this case. Consequently, the whole proof of the lemma is completed.
\qed

		Next, we prove a continuation lemma based on the a priori estimate in Lemma \ref{lem:3.2}.
		\begin{lem}\label{lem:3.3}

			Under the condition of Assumption  \ref{ass:3.1} and Assumption \ref{ass:3.2}, if for some $\alpha_{0}\in[0,1)$, FBS$\bigtriangleup $E \eqref{eq:3.8} admits a unique solution $(x(\cdot),y(\cdot))\in N^2(\overline{\mathbb{T}};\mathbb{R}^{2n} )$ for any $(\xi,\eta ,\rho (\cdot ))\in \mathcal{H} (\overline{\mathbb{T}})$, then there exists an absolute constant $\delta_0>0$ such that FBS$\bigtriangleup $E \eqref{eq:3.8} admits a unique solution   for $\alpha=\alpha_{0}+\delta$ with $\delta\in(0,\delta_0]$ and $\alpha\le1$.
		\end{lem}
		$\mathbf{Proof.}$
			Let $\delta_0>0$ be determined below. For any $(x(\cdot),y(\cdot))\in N^2(\overline{\mathbb{T}};\mathbb{R}^{2n} )$, $\theta(k)=(x_k,y'_{k+1},z'_{k+1}) $ and $ \Theta(k)=(X_k,Y'_{k+1},Z'_{k+1})$,we introduce the following FBS$\bigtriangleup $E with unknow $(X(\cdot),Y(\cdot))\in N^2(\overline{\mathbb{T}};\mathbb{R}^{2n} )$ :
			\begin{equation}\label{eq:3.25}
				\left\{\begin{aligned}
					X_{k+1}  = &\Big\{-(1-\alpha_0)\mu\Big[B_k^{\top}Q(k,  Y_{k+1}', Z'_{k+1})\Big] \\
					& \quad + \alpha_0 b(k, \Theta(k)) + \widetilde{\psi}_k\Big\} \\
					& \quad + \Big\{-(1-\alpha_0)\mu\Big[C_k^{\top}Q(k, Y_{k+1}', Z'_{k+1})\Big] \\
					& \quad + \alpha_0 \sigma(k, \Theta(k)) + \widetilde{\gamma}_k\Big\} w_k, \\
					Y_k &= -\Big\{-(1-\alpha_0)v\Big[A_k^{\top}P(k, X_k)\Big] \\
					& \quad + \alpha_0 f(k, \Theta(k)) + \widetilde{\varphi}_k\Big\}, \quad k \in \mathbb{T}, \\
					X_0 &= \alpha_{0}\Lambda(Y_0) - (1-\alpha_{0})\mu M^{\top}MY_0 + \widetilde{\xi}, \\
					Y_N &= \alpha_{0}\Phi(X_N) + (1-\alpha_{0})vG^{\top}GX_N + \widetilde{\eta},
				\end{aligned}\right.
			\end{equation}
			where
			\begin{equation}
				\left\{\begin{aligned}
					&\widetilde{\psi}_k=\psi_k +\delta b(k,\theta (k))+\delta \mu  \Big [B_k^\top Q(k, y_{k+1}', z'_{k+1})\Big ],\\
					&\widetilde{\gamma}_k=\gamma_k +\delta \sigma(k,\theta (k))+\delta \mu  \Big [C_k^\top Q(k, y_{k+1}', z'_{k+1}))\Big ],\\
					&\widetilde{\varphi}_k=\varphi_k+\delta f(k+1,\theta (k))+\delta v\Big[A_k^\top P(k,x_k)\Big],\\
					&\widetilde{\eta}=\eta +\delta \Phi(x_N)-\delta vG^\top Gx_N,\\
					&\widetilde{\xi}=\xi+\delta\Lambda(y_0)+\delta\mu M^\top My_0,\\
					&Z_{k+1}=Y_{k+1}\omega_k,\\
					&Y'_{k+1}=[Y_{k+1}|\mathcal{F}_{k-1}],\\
					&Z'_{k+1}=[Y_{k+1}\omega_k|\mathcal{F}_{k-1}].
				\end{aligned}\right.
			\end{equation}
			Furthermore, we also denote   $\widetilde{\rho} (\cdot )=(\widetilde{\varphi} (\cdot ) ,\widetilde{\psi} (\cdot ),\widetilde{\gamma} (\cdot ))$. Then it is easy to check that $(\widetilde{\xi},\widetilde{\eta},\widetilde{\rho})\in\mathcal{H}[\overline{\mathbb{T}}]$. By our assumptions, the FBS$\bigtriangleup $E \eqref{eq:3.25} admits a unique solution $(X(\cdot),Y(\cdot))\in N^2(\overline{\mathbb{T}};\mathbb{R}^{2n} )$. In fact, we can  establish a mapping
			\begin{equation}
				\begin{aligned}
					&(X(\cdot),Y(\cdot))
					\\&=\mathcal{T}_{\alpha _0+\delta }\big ((x(\cdot),y(\cdot))\big ): N^2(\overline{\mathbb{T}};\mathbb{R}^{2n} )\to N^2(\overline{\mathbb{T}};\mathbb{R}^{2n} ).\nonumber
				\end{aligned}
			\end{equation}
			In the following, we shall prove that the above mapping is contractive when $\delta$ is small enough.

			Let $(x(\cdot),y(\cdot)),(\bar{x}(\cdot),\bar{y}(\cdot))\in N^2(\overline{\mathbb{T}};\mathbb{R}^{2n} )$ and $(X(\cdot),Y(\cdot))=\mathcal{T}_{\alpha _0+\delta }\big (x(\cdot),y(\cdot)\big )$, $(\bar{X}(\cdot),\bar{Y}(\cdot))=\mathcal{T}_{\alpha _0+\delta }\big (\bar{x} (\cdot ),\bar{y}(\cdot)\big )$. Similarly, we denote $\widehat{x }(\cdot)=x(\cdot)-\bar{x}(\cdot),\widehat{y }(\cdot)=y(\cdot)-\bar{y}(\cdot),\widehat{X }(\cdot)=X(\cdot)-\bar{X}(\cdot),\widehat{Y }(\cdot)=Y(\cdot)-\bar{Y}(\cdot)$, etc. By applying Lemma \ref{lem:3.2}, we have (the argument $k$ is suppressed for simplicity)
			\begin{equation}
				\begin{aligned}
					&\quad\mathbb E\bigg[\displaystyle \sum_{k=0}^N |\widehat{X}_k|^2+\displaystyle \sum_{k=0}^N |\widehat{Y}_k|^2\bigg ]\\
					\le & K\delta ^2\mathbb{E}\Bigg \{\Big |\big (\Phi (x_N)-\Phi (\bar{x}_N)\big )-v G^\top G\widehat{x}_N\Big |^2\\
					&\qquad \quad +\Big|\Lambda(y_0)-\Lambda(\bar{y}_0)-\mu M^\top M\widehat{y_0}    \Big|^2\\
					&\quad +\sum_{k=0}^{N-1}\bigg |\big (b(\theta_k)-b(\bar{ \theta_k })\big )+\mu \big[B_k^\top Q(k,\widehat y_{k+1}',\widehat z_{k+1})\big]\bigg |^2\\
					&\quad +\sum_{k=0}^{N-1}\bigg |\big (\sigma(\theta_k)-\sigma(\bar{ \theta_k }\big )+\mu \big[C_k^\top Q(k,\widehat y_{k+1}',\widehat z_{k+1})\big]\bigg |^2\\
					&\quad+\sum_{k=0}^{N-1}\bigg|\big (f(\theta_k)-f(\bar{ \theta_k })\big )+v\big[A_k^\top P(k,\widehat{x}_k)\big] \bigg|
					\Bigg\}.\nonumber
				\end{aligned}
			\end{equation}
			Due to the Lipschitz continuity of $(\Lambda,\Phi,\Gamma)$ and the boundedness of $G,M,A_k,B_k,C_k$, there exists a new constant $K'>0$ independent of $\alpha_{0}$ and $\delta$ such that
			\begin{equation}
				\|(\widehat{X}(\cdot),\widehat{Y}(\cdot)) \|^2_{N^2(\overline{\mathbb{T}};\mathbb{R}^{2n} )}\le K'\delta ^2 \|(\hat{x}(\cdot),\hat{y}(\cdot)) \|^2_{N^2(\overline{\mathbb{T}};\mathbb{R}^{2n} )}.\nonumber
			\end{equation}
			By selecting $\delta_0=1/(2\sqrt{K'})$, when $\delta\in(0,\delta_0]$, we can find that the mapping $\mathcal{T}_{\alpha _0+\delta }$ is contractive. Thus, the mapping $\mathcal{T}_{\alpha _0+\delta }$ admits a unique fixed point which is just the unique solution to FBS$\bigtriangleup $E \eqref{eq:3.8}. The proof is finished.
\qed

		We summarize the above analysis to give the following proof.  \qed

		$\mathbf{[Proof~of~Theorem ~ \ref{thm:3.1}.]}$
			Firstly, the  solvability of FBS$\bigtriangleup $E \eqref{eq:1.1} in the space $N^2(\overline{\mathbb{T}};\mathbb{R}^{2n} )$ is deduced from the unique solvability of FBS$\bigtriangleup $E \eqref{eq:3.10} and Lemma \ref{lem:3.3}. Secondly, in Lemma \ref{lem:3.2}, by taking $\alpha=1$, $\big(\xi,\eta,\rho(\cdot)\big)=(0,0,0)$, and$\big(\bar{\xi},\bar{\eta},\bar{\rho}(\cdot)\big)=\Big(\bar{\Lambda}(\bar{y}_0)-\Lambda(\bar{y}_0), \bar{\Phi}(\bar{x}_N)-\Phi(\bar{x}_N), \bar{\Gamma}\big(\cdot,\bar{\theta}(\cdot)\big)-\Gamma\big(\cdot,\bar{\theta}(\cdot)\big)\Big)$, we get the estimate \eqref{eq:3.6} in Theorem \ref{thm:3.1} from the estimate \eqref{eq:3.11} in Lemma \ref{lem:3.2}. Finally, by selecting the coefficients $\big(\bar{\Lambda},\bar{\Phi},\bar{\Gamma}\big)=(0,0,0)$, we get \eqref{eq:3.4} from \eqref{eq:3.6}. For the uniqueness of the solution of FBS$\bigtriangleup $E(1.1) , it can be obtained directly from the a priori estimate (\ref{eq:3.6}). The proof is completed. \qed

 To illustrate the applicability of the results, two linear-quadratic (LQ) examples are presented in the next section.

	\section{Applications in linear quadratic problem}\label{sec:4}

		In this section, we aim to establish the existence and uniqueness of solutions for stochastic Hamiltonian systems associated with certain types of linear-quadratic (LQ) optimal control problems. If it can be shown that these stochastic Hamiltonian systems adhere to the previously introduced domination-monotonicity conditions, then, according to Theorem \ref{thm:3.1}, we can conclude that their solutions exist and are unique. The initial motivation behind this paper's investigation into fully coupled FBS$\bigtriangleup $Es (\ref{eq:1.1}) is to determine the solvability of these Hamiltonian systems.

		\subsection{Forward LQ stochastic control problem}\label{sec:4.1}

		In this section, we consider the following linear forward S$\bigtriangleup $E control system:

		\begin{equation}\label{eq:4.1}
			\left\{\begin{aligned}
				x_{k+1}&=A_k x_k+B_k u_k+b_k+\left(C_k x_k+D_k u_k+\sigma_k\right) \omega_k, \\
				x_0&=\xi ,\quad k \in \mathbb{T}
			\end{aligned}\right.
		\end{equation}
		where $b(\cdot),\sigma(\cdot)\in L_{\mathbb{F}}^2\left(\mathbb{T}; \mathbb{R}^{n}\right);A(\cdot),C(\cdot) \in L_{\mathbb{F}}^{\infty}(\mathbb{T}; \mathbb{R}^{n\times n}),\\ B(\cdot),D(\cdot) \in L_{\mathbb{F}}^{\infty}(\mathbb{T}; \mathbb{R}^{n\times m}). $  The vector $\xi \in \mathbb{R}^n$ is called an initial state.  Here the stochastic process  $u(\cdot)=\left(u_{0}, u_{1}, \ldots, u_{N-1}\right)  \in \mathbb{U}$ is called  an admissible control. By Lemma \ref{lem:2.2}, we know that  the S$\bigtriangleup $E (\ref{eq:4.1}) admits a unique solution $x(\cdot )\in L_{\mathbb{F}}^2\left(\overline{\mathbb{T}}; \mathbb{R}^n\right)$.

		Meanwhile, we propose a quadratic criterion function:
		\begin{equation}\label{eq:4.2}
			\begin{aligned}
				J\left(0; \xi, u(\cdot)\right) &= \frac{1}{2} \mathbb{E} \Bigg\{ \left\langle M x_0, x_0 \right\rangle + \left\langle G x_N, x_N \right\rangle \\
				& \quad + \sum_{k=0}^{N-1} \left[ \left\langle Q_k x_k, x_k \right\rangle + \left\langle R_k u_k, u_k \right\rangle \right] \Bigg\}.
			\end{aligned}
		\end{equation}
		where $M\in \mathbb{S}^n , G\in L^2 _{\mathcal{F}_{N-1} }(\Omega ;\mathbb{S}^n),
		Q(\cdot) \in  L_{\mathbb{F}}^{\infty}(\mathbb{T}; \mathbb{S}^{n}) $ and  $R(\cdot) \in L_{\mathbb{F}}^{\infty}(\mathbb{T}; \mathbb{S}^{m})$.

		Next, we give the main questions to be explored as follows:

		\textbf{Problem(FLQ).} The problem is to find a pair of initial state and admissible control $(\bar \xi, \bar u(\cdot)) \in  \mathbb{R}^n \times \mathbb{U} $ such that

		\begin{equation}
			J(0; \bar\xi,  \bar{u}(\cdot))=\inf _{(\xi, u(\cdot)) \in \mathbb{R}^n \times \mathbb{U}} J\left(0, \xi, u(\cdot)\right).
		\end{equation}
		Problem(FLQ) is said to be solvable if there exists $(\bar\xi,  \bar{u}(\cdot))$ above such that infimum of ceriterion function is achieved. Then $\bar{u}(\cdot)$ is called an optimal conrol and $\bar{\xi}$ is called an optimal initial state.

		Moreover, we impose the following assumption:

		\begin{ass}\label{ass:4.1}
			(i)$M$ is positive definite and  for any $\omega \in \Omega$ , $G(\omega)$ is nonnegative definite;\\
			(ii)For any $(\omega ,k)\in \Omega \times  \mathbb{T}$, $Q_k$ is nonnegative definite;\\
			(iii)For any $(\omega ,k)\in \Omega \times  \mathbb{T}$, there exixts a    constant $\delta >$ 0 and unit matrix $I_m$ such that  $R_k - \delta I_m \ge 0$ .
		\end{ass}

		After the above description of the problem and the formulation of assumptions,  we can use Theorem \ref{thm:3.1} to study Problem(FLQ). For  admissible control $(\bar \xi, \bar u(\cdot)) \in  \mathbb{R}^n \times \mathbb{U} ,$ we introduce the stochastic Hamiltonian system of S$\bigtriangleup $E(\ref{eq:4.1}) as follows:

		\begin{equation} \label{eq:4.6}
			\left\{\begin{aligned}
				\bar{x}_{k+1}=&A_k \bar{x}_k+B_k \bar{u}_k+b_k+\left\{C_k \bar{x}_k+D_k \bar{u}_k+\sigma_k\right\} \omega_k, \\
				\bar{y}_k=&A_k^{\top} \bar{y}_{k+1}^{\prime}+C_k^{\top} \bar{z}_{k+1}^{\prime}+Q_k \bar{x}_k, \quad k \in \mathbb{T}, \\
				\bar{x}_0=&-M^{-1} \bar{y}_0  ,\\
				\bar{y}_N=&G \bar{x}_N ,\\
				0=&B_k^{\top} \bar{y}_{k+1}^{\prime}+D_k^{\top} \bar{z}_{k+1}^{\prime}+ R_k\bar{u}_k .
			\end{aligned}\right.
		\end{equation}

		It is clear that
		the Hamiltonian system \eqref{eq:4.6}
		is described by an FBS$\bigtriangleup $E satisfying Assumptions \ref{ass:3.1} and \ref{ass:3.2}. Therefore, with the help of Theorem \ref{thm:3.1}, we can get the following theorem.

		\begin{thm}\label{thm:4.1}
			Under Assumption \ref{ass:4.1}, the above Hamiltonian system \eqref{eq:4.6} admits a unique solution $(x(\cdot ),y(\cdot )) \in N^{2}(\overline{\mathbb{T}} ;\mathbb{R}^{2n} )$. Moreover,\\

			\begin{equation}\label{eq:4.7}
				\bar{\xi}=-M^{-1} \bar{y}_0
			\end{equation}
			is the unique optimal initial state and

			\begin{equation}\label{eq:4.8}
				\bar{u}_k=-R_k^{-1}\left(B_k^{\top} \bar{y}_{k+1}^{\prime}+D_k^{\top} \bar{z}_{k+1}^{\prime}\right)
			\end{equation}
			is the unique optimal control of Problem(FLQ)
		\end{thm}
		$\mathbf{{Proof.}}$ Firstly, since  $R_k$ is invertible,  we can solve the last equation of Hamiltonian system (\ref{eq:4.6}) to get
			\begin{equation} \label{eq:4.9}
				\bar{u}_k=-R_k^{-1}\left(B_k^{\top} \bar{y}_{k+1}^{\prime}+D_k^{\top} \bar{z}_{k+1}^{\prime}\right) .
			\end{equation} \\
			By substituting equation (\ref{eq:4.9}) into the Hamiltonian system (\ref{eq:4.6}), we derive a forward-backward stochastic difference equation (FBS$\bigtriangleup $E). We then verify that the coefficients of this FBS$\bigtriangleup $E satisfy both Assumptions \ref{ass:3.1} and \ref{ass:3.2}. Consequently, Theorem \ref{thm:3.1} guarantees the existence of a unique solution to this FBS$\bigtriangleup $E, which in turn ensures a unique solution to the Hamiltonian system (\ref{eq:4.6}).

			Consider $(\bar \xi, \bar u(\cdot))\in \mathbb{R}^n \times \mathbb{U}$, which satisfies conditions (\ref{eq:4.7}) and (\ref{eq:4.8}). We will now prove the optimality and uniqueness of the pair (\ref{eq:4.7}) and (\ref{eq:4.8}). We begin by establishing the existence of an optimal initial state and optimal control. Let $(\xi,u(\cdot)) \in \mathbb{R}^n \times \mathbb{U}$ be an arbitrary set of initial states and admissible controls, distinct from $(\bar \xi, \bar u(\cdot))$, and let $x(\cdot)$ be the corresponding state process for $(\xi,u(\cdot))$. We will now analyze the difference in the criterion function between these two distinct states:
		{\scriptsize	\begin{equation} \label{eq:4.10}
				\begin{aligned}
					& J\left(0 ; \xi, u_k\right)-J(0;\bar{\xi}, \bar{u}_k) \\
					=&\frac{1}{2} \mathbb{E} \Bigg \{\left\langle M \xi, \xi \right\rangle-\langle M \bar{\xi}, \bar{\xi}\rangle+\left\langle G x_N, x_N\right\rangle-\left\langle G \bar{x}_N, \bar{x}_N\right\rangle \\
									&\left.+\sum_{k=0}^{N-1}\bigg[\left\langle Q_k x_k, x_k\right\rangle-\left\langle Q_k \bar{x}_k, \bar{x}_k\right\rangle+\left\langle R_k u_k, u_k\right\rangle-\left\langle R_k \bar{u}_k, \bar{u}_k\right\rangle\bigg]\right\} \\
					=&\frac{1}{2} \mathbb{E} \Bigg \{\left\langle M\left(\xi-\bar{\xi}) \right),\xi-\bar{\xi} \right\rangle+\left\langle G\left(x_N-\bar{x}_N\right), x_N-\bar{x}_N\right\rangle \\
					&+\sum_{k=0}^{N-1}\bigg[\left\langle Q_k\left(x_k-\bar{x}_k\right), x_k-\bar{x}_k\right\rangle+\left\langle R_k\left(u_k-\bar{u}_k\right), u_k-\bar{u}_k\right\rangle\bigg]\bigg\}\\&+\Delta_1,
				\end{aligned}
			\end{equation}}
			where
			\begin{equation}\label{eq:901}
				\begin{aligned}
					\Delta_1=&\mathbb{E} \Bigg \{  \left\langle M \bar{\xi}, \xi-\bar{\xi} \right\rangle+\left\langle G \bar{x}_N, x_N-\bar{x}_N\right\rangle
					\\& +\sum_{k=0}^{N-1}\left[\left\langle Q_k \bar{x}_k, x_k-\bar{x}_k\right\rangle+\left\langle R_k \bar{u}_k, u_k-\bar{u}_k\right\rangle\right] \Bigg \}.
				\end{aligned}
			\end{equation}
			Furthermore,   based on the expressions for the initial value and terminal value in \eqref{eq:4.6}, we obtain
			\begin{equation} \label{eq:4.101}
				\begin{aligned}
					&\mathbb{E} \Bigg \{  \left\langle M \bar{\xi}, \xi-\bar{\xi} \right\rangle+\left\langle G \bar{x}_N, x_N-\bar{x}_N\right\rangle\bigg\}
					\\=&\mathbb{E} \Bigg \{\left\langle x_N-\bar{x}_N, \bar{y}_N\right\rangle-\left\langle \xi-\bar{\xi}, \bar{y}_0\right\rangle \Bigg \} \\
					=&   \mathbb{E} \Bigg \{\sum_{k=0}^{N-1}\bigg[\left\langle x_{k+1}-\bar{x}_{k+1}, \bar{y}_{k+1}\right\rangle-\left\langle x_k-\bar{x}_k, \bar{y}_k\right\rangle   \bigg]\bigg\}
					\\= &\mathbb{E}\bigg\{\sum_{k=0}^{N-1}\bigg[\Big\langle A_k(x_k -\bar{x}_k )+B_k(u_k-\bar{u}_k)\\
					&\qquad +\big[ C_k(x_k -\bar{x}_k)+D_k(u_k-\bar{u}_k)\big] \omega_k, \bar{y}_{k+1}\Big\rangle \\
					& \qquad -\left\langle x_k-\bar{x}, A_k^{\top} \bar{y}_{k+1}^{\prime}+C_k^{\top} \bar{z}_{k+1}^{\prime}+Q_k \bar{x}_k\right\rangle \bigg]\Bigg \} \\
					=&\mathbb{E} \Bigg \{ \sum_{k=0}^{N-1}\bigg[\left\langle B_k^{\top} \bar{y}_{k+1}^{\prime}+D_k^{\top} \bar{z}_{k+1}^{\prime}, u_k-\bar{u}_k\right\rangle\\
					&\qquad \quad -\left\langle Q_k \bar{x}_k , x_k-\bar{x}_k\right\rangle\bigg]\Bigg \}.
			\end{aligned}\end{equation}
			Putting \eqref{eq:4.101} into \eqref{eq:901}, we get
			\begin{equation}\label{eq:411}
				\begin{aligned}
					\Delta_1=& \mathbb{E} \Bigg \{ \sum_{k=0}^{N-1}\left\langle B_k^{\top} \bar{y}_{k+1}^{\prime}+D_k^{\top} \bar{z}_{k+1}^{\prime}+R_k \bar{u}_k, u_k-\bar{u}_k\right\rangle \Bigg \}
					\\=&0 ,
				\end{aligned}
			\end{equation}
			where we have used  $$0=B_k^{\top} \bar{y}_{k+1}^{\prime}+D_k^{\top} \bar{z}_{k+1}^{\prime}+ R_k\bar{u}_k $$ in  \eqref{eq:4.6}.
			Then putting \eqref{eq:411} into \eqref{eq:4.10} and by the positive definiteness of the relevant coefficient matrices in Assumption \ref{ass:4.1}, we can easily verify that $$J\left(0 ; \xi, u(\cdot)\right)-J\left(0 ; \bar{\xi}, \bar{u}(\cdot)\right) \ge 0.$$  Since ($\xi ,u(\cdot)$) is an arbitrary pair of initial state and admissible control, we obtain the optimality of ($\bar{\xi},\bar{u}(\cdot)$).

			Next, we demonstrate the uniqueness.
			
			Let $(\tilde{x}_0, \tilde{v}(\cdot))\in \mathbb{R}^n \times \mathbb{U}$ represent another optimal initial state and control pair, apart from $(\bar{\xi},\bar{u}(\cdot))$. The corresponding optimal state process for $(\tilde{x}_0, \tilde{v}(\cdot))$ is denoted by $\tilde{x}(\cdot)$. This implies the condition $J(0 ; \bar{\xi}, \bar{u}(\cdot)) = J(0 ; \tilde{x}_0, \tilde{v}(\cdot))$. Referring back to (\ref{eq:4.10}), we can infer:
			\begin{equation}
				\begin{aligned}
					0 = &\ J(0; \bar{\xi}, \bar{u}(\cdot)) - J(0; \tilde{x}_0, \tilde{v}(\cdot)) \\
					= &\ \frac{1}{2} \mathbb{E} \Bigg\{ \left\langle M\left(\bar{\xi} - \tilde{x}_0\right), \bar{\xi} - \tilde{x}_0 \right\rangle \\
					&\qquad + \left\langle G\left(\bar{x}_N - \tilde{x}_N\right), \bar{x}_N - \tilde{x}_N \right\rangle \\
					&\qquad + \sum_{k=0}^{N-1} \bigg[ \left\langle Q_k\left(\bar{x}_k - \tilde{x}_k\right), \bar{x}_k - \tilde{x}_k \right\rangle \\
					&\qquad+ \left\langle R_k\left(\bar{u}_k - \tilde{v}_k\right), \bar{u}_k - \tilde{v}_k \right\rangle \bigg] \Bigg\} \\
					\geq &\ \frac{1}{2} \mathbb{E} \Bigg\{ \left\langle M\left(\bar{\xi} - \tilde{x}_0\right), \bar{\xi} - \tilde{x}_0 \right\rangle \\
					&\qquad + \sum_{k=0}^{N-1} \left\langle R_k\left(\bar{u}_k - \tilde{v}_k\right), \bar{u}_k - \tilde{v}_k \right\rangle \Bigg\}.
				\end{aligned}
			\end{equation}
			Due to positive definiteness of $M$ and $R_k$, we have $\bar{\xi} = \tilde{x}_0$ and $\bar{u}_k=\tilde{v}_k$. The uniqueness is proved. The proof is completed.
\qed

		\subsection{Backward LQ stochastic control problem}\label{sec:4.2}

		In this section, we consider the following linear BS$\bigtriangleup $E control system:
		\begin{equation}
			\left\{\begin{aligned} \label{eq:4.14}
				y_k=&A_k y'_{k+1}+B_k z'_{k+1}+C_k v_k+\alpha_k,  \quad k\in \mathbb{T},\\
				y_N=&\eta,
			\end{aligned}\right.
		\end{equation}
		where $A(\cdot),B(\cdot)\in L_{\mathbb{F}}^{\infty}(\mathbb{T}; \mathbb{R}^{n\times n}),
		C(\cdot) \in L_{\mathbb{F}}^{\infty}(\mathbb{T}; \mathbb{R}^{n\times m})$,\\ $\alpha(\cdot) \in L_{\mathbb{F}}^2\left(\mathbb{T}; \mathbb{R}^{n}\right). $   $\eta \in L^{2} _{\mathcal{F}_{N-1} }(\Omega ; \mathbb {R}^n)$ is called  terminal state.  The process $v_k$ is called admissible control and $v(\cdot)=\left(v_{0}, v_{1}, \ldots, v_{N-1}\right)  \in \mathbb{U}$. By Lemma \ref{lem:2.3}, we know that the solution of the BS$\bigtriangleup $E (\ref{eq:4.14}) admits a unique solution $y(\cdot )\in L_{\mathbb{F}}^2\left(\overline{\mathbb{T}}; \mathbb{R}^n\right)$.

		Meanwhile, we propose a quadratic criterion function:
		\begin{equation}  \label{eq:4.131}
			\begin{aligned}
				& J(N ; \eta, v(\cdot)) \\
				&=\frac{1}{2} \mathbb{E}\bigg\{\left\langle M y_0, y_0\right\rangle+\sum_{k=0}^{N-1}\bigg[\left\langle Q_{k} y'_{k+1} ,y'_{k+1}\right\rangle
\\
				&\qquad \quad +\left\langle L_k z'_{k+1}, z'_{k+1}\right\rangle+\left\langle R_k v_k, v_k\right\rangle \bigg]\bigg\},
			\end{aligned}
		\end{equation}
		where $M\in \mathbb{S}^n ,
		Q(\cdot),L(\cdot) \in    L_{\mathbb{F}}^{\infty}(\mathbb{T}; \mathbb{S}^n) $ and  $R(\cdot) \in L_{\mathbb{F}}^{\infty}(\mathbb{T}; \mathbb{S}^m) $.

		Next, we propose the following problem:

		\textbf{Problem(BLQ).} For any given terminal state $\eta \in L^{2} _{\mathcal{F}_{N-1} }(\Omega ; \mathbb {R}^n)$, find an admissible control $\bar{v}(\cdot) \in \mathbb{U}$ such that

		\begin{equation}
			J(N ; \eta, \bar{v}(\cdot))=\inf _{v(\cdot) \in M^2\left(T ; R^m\right)} J(N ; \eta, v(\cdot)).
		\end{equation}

		Problem(BLQ) is said to be solvable if there exists $\bar{v}(\cdot)$ above such that the infimum of the ceriterion function \eqref{eq:4.131} is achieved. Then $\bar{v}(\cdot)$ is called an optimal control.

		Moreover, we impose the following assumption:
		\begin{ass}\label{ass:4.2}
			(i)$M$ is positive ;\\
			(ii)For any $(\omega ,k)\in \Omega \times  \mathbb{T}$, $Q(\cdot)$ and $L(\cdot)$ are nonnegative definite;\\
			(iii)For any $(\omega ,k)\in \Omega \times  \mathbb{T}$, there exixts a    constant $\delta >$ 0 and unit matrix $I_m$ such that  $R_k - \delta I_m \ge 0$ .
		\end{ass}
		Hamiltonian system for Problem(BLQ) given some possible admissible
		control pair $(\bar v(\cdot), \bar x(\cdot))$ as follows:
		\begin{equation}\label{eq:4.17}
			\left\{\begin{aligned}
				\bar x_{k+1}&=A_k^{\top} \bar x_k-Q_k \bar y'_{k+1}+\left(B_k^{\top}\bar  x_k-L_k \bar z'_{k+1}\right) \omega_k , \quad k\in
				\mathbb{T} ,\\
				\bar y_k&=A_k \bar y'_{k+1}+B_k \bar z'_{k+1}+C_k \bar v_k+\alpha_k , \\
				\bar x_0&=-M \bar y_0 ,\\
				\bar y_N&=\eta ,\\
				0&=C_k^{\top} \bar x_k-R_k \bar v_k.
			\end{aligned}\right.
		\end{equation}
		\begin{thm}\label{thm:4.2}
			Under Assumption \ref{ass:4.2}, the above Hamiltonian system \eqref{eq:4.17} admits a unique solution $(\bar x(\cdot ),\bar y(\cdot ) ) \in N^{2}(\overline{\mathbb{T}} ;\mathbb{R}^{2n} )$. Moreover, for a given terminal state $\eta \in L^{2}_{\mathcal{F}_{N-1}}(\Omega ; \mathbb{R}^n)$ ,
			\\
			\begin{equation} \label{eq:4.18}
				\bar v_k=R^{-1}_k C_k^{\top} \bar x_k
			\end{equation}
			is the unique optimal control of Problem(BLQ).
		\end{thm}
		$\mathbf{{Proof.}}$
			From the final relationship in \eqref{eq:4.17}, we obtain the relationship (\ref{eq:4.18}). Substituting (\ref{eq:4.18}) into the Hamiltonian system (\ref{eq:4.17}), we derive an FBS$\bigtriangleup $E. We can then verify that the coefficients of this FBS$\bigtriangleup $E satisfy Assumption \ref{ass:3.1} and Assumption \ref{ass:3.2}. Therefore, by Theorem \ref{thm:3.1}, there exists a unique solution $(\bar x(\cdot ),\bar y(\cdot ) ) \in N^{2}(\overline{\mathbb{T}} ;\mathbb{R}^{2n} )$ to this FBS$\bigtriangleup $E, implying a unique solution to the Hamiltonian system (\ref{eq:4.17}).

			Now we begin to prove  (\ref{eq:4.18}) is a unique optimal control.
			Let $v(\cdot) \in  \mathbb{U} $ be any given  admissible control except $\bar{v}(\cdot) $ and $ y(\cdot) $ be the corresponding state process with $v(\cdot)$. We consider the difference between the  cost functional in these two different states:
			\begin{equation} \label{eq:4.19}
				\begin{aligned}
					& J(N ; \eta, v(\cdot))-J(N ; \eta, \bar{v}(\cdot)) \\
					& =\frac{1}{2} \mathbb{E}\Bigg\{ \left\langle M y_0, y_0\right\rangle-\left\langle M \bar{y}_0, \bar{y}_0\right\rangle\\
										&\qquad \quad +\sum_{k=0}^{N-1}\bigg[\left\langle Q_k {y}'_{k+1}, {y}'_{k+1}\right\rangle -\left\langle Q_k \bar{y}'_{k+1}, \bar{y}'_{k+1}\right\rangle\\
										& \left. \qquad \qquad \qquad +\left\langle L_k {z}'_{k+1}, {z}'_{k+1}\right\rangle -\left\langle L_k \bar{z}'_{k+1}, \bar{z}'_{k+1}\right\rangle\right.\\
									&\qquad \qquad \qquad +\left. \left\langle R_k u_k, u_k\right\rangle -\left\langle R_k \bar{u}_k, \bar{u}_k\right\rangle\bigg]\right\} \\
					& =\frac{1}{2}\mathbb{E}\Bigg\{\left\langle M\left(y_0-\bar{y}_0\right), y_0-\bar{y}_0\right\rangle \\
					&\qquad+\sum_{k=0}^{N-1}\bigg[\left\langle Q_k\left({y}'_{k+1}-\bar{y}'_{k+1}\right), {y}'_{k+1}-\bar{y}'_{k+1}\right\rangle\\
					&\qquad+\left\langle L_k\left({z}'_{k+1}-\bar{z}'_{k+1}\right), {z}'_{k+1}-\bar{z}'_{k+1}\right\rangle
					\\&\qquad+ \left\langle R_k(v_k-\bar{v}_k),v_k-\bar{v}_k\right\rangle \bigg]\bigg\}+\Delta_2,
\end{aligned}
			\end{equation}
where
\begin{equation}
				\begin{aligned}
					\Delta_2=\mathbb{E} \bigg\{&\langle M \bar{y}_0, y_0-\bar{y}_0\rangle\\
					&+\sum_{k=0}^{N-1}\bigg[\langle Q_k {\bar{y}}'_{k+1}, {y}'_{k+1}-\bar{y}'_{k+1}\rangle
					\\&+\langle L_k {\bar{z}}'_{k+1},  {z}'_{k+1}-\bar{z}'_{k+1} \rangle\\
					&\qquad +\langle R_k \bar{v}_k, v_k-\bar{v}_k\rangle\bigg]\bigg\}.
				\end{aligned}
			\end{equation}
			Furthermore,  from the relationships in (\ref{eq:4.17}) and noting that $y_N-\bar{y}_N=\eta-\eta=0$ , we can get that
			\begin{equation*} \label{eq:4.191}
\begin{split}
&\mathbb{E}[\left\langle M \bar{y}_0, y_0-\bar{y}_0\right\rangle] \\
&= \mathbb E\bigg\{\sum_{k=0}^{N -1}\bigg[\left\langle \bar{x}_{k+1}, y_{k+1}-\bar{y}_{k+1}\right\rangle-\left\langle \bar{x}_k, y_{k}-\bar{y}_{k}\right\rangle\bigg]\bigg\} \\
&= \mathbb{E}\bigg\{ \sum_{k=0}^{N-1}\bigg[
\langle A_k^{\top} \bar x_k-Q_k \bar y'_{k+1}+(B_k^{\top}\bar x_k-L_k \bar z'_{k+1}) \omega_k, \\
&\qquad y_{k+1}-\bar{y}_{k+1}\rangle
- \langle \bar{x}_k , A_k (y'_{k+1}-\bar y'_{k+1}) \\
&\quad +B_k (z'_{k+1}-\bar z'_{k+1})+C_k (v_k-\bar v_k)\rangle
\bigg]\bigg\} \\
&= \mathbb{E}\bigg\{\sum_{k=0}^{N-1}\bigg[-\langle Q_k \bar{y}'_{k+1}, y'_{k+1}-\bar{y}'_{k+1}\rangle \\
&\quad -\langle L_k \bar{z}'_{k+1}, z'_{k+1}-\bar{z}'_{k+1}\rangle-\langle C_k^{\top} \bar{x}_k, v_k-\bar{v}_k\rangle\bigg] \bigg\} \\
&= \mathbb{E}\bigg\{\sum_{k=0}^{N-1}\bigg[-\langle Q_k \bar{y}'_{k+1}, y'_{k+1}-\bar{y}'_{k+1}\rangle
\end{split}
\end{equation*}
\begin{equation} \label{eq:4.191}
\begin{split}
&\quad -\langle L_k \bar{z}'_{k+1}, z'_{k+1}-\bar{z}'_{k+1}\rangle
- \langle R_k^{\top} \bar{v}_k, v_k-\bar{v}_k\rangle\bigg]\bigg\}.
\end{split}
\end{equation}
			As a result, $\Delta_2 = 0$. Then, from \eqref{eq:4.19},  taking into account the positive definiteness of the relevant coefficient matrices in Assumption \ref{ass:4.2}, it can be easily shown that $J(N ; \eta, v(\cdot)) - J(N ; \eta, \bar{v}(\cdot)) \ge 0$. Because $v(\cdot)$ is an arbitrary admissible control, this establishes the optimality of $\bar{v}(\cdot)$.

			We now proceed to prove uniqueness. For the given terminal state $\eta \in L^{2}_{\mathcal{F}_{N-1}}(\Omega ; \mathbb{R}^n)$, assume that there exists another optimal control $\tilde{u}(\cdot) \in \mathbb{U}$ in addition to $\bar{v}(\cdot)$. Let $\tilde{y}(\cdot)$ denote the corresponding state process for $\tilde{u}(\cdot)$. This implies that $J(N ; \eta, \tilde{u}(\cdot)) = J(N ; \eta, \bar{v}(\cdot))$. Then, referring back to (\ref{eq:4.19}) and \eqref{eq:4.191}, we can deduce the following:
			{\small	\begin{equation}
					\begin{aligned}
						0= &J(N ; \eta, \tilde{u}(\cdot)) - J(N ; \eta, \bar{v}(\cdot))
						\\=&\frac{1}{2}\mathbb{E}\Bigg\{\left\langle M\left(\tilde{y}_0-\bar{y}_0\right), \tilde{y}_0-\bar{y}_0\right\rangle
						\\&\qquad+\sum_{k=0}^{N-1}\bigg[\left\langle Q_k\left(\tilde{y}'_{k+1}-\bar{y}'_{k+1}\right), \tilde{y}'_{k+1}-\bar{y}'_{k+1}\right\rangle\\
						&\qquad+\left\langle L_k\left(\tilde{z}'_{k+1}-\bar{z}'_{k+1}\right), \tilde{z}'_{k+1}-\bar{z}'_{k+1}\right\rangle
						\\&\qquad+ \left\langle R_k(\tilde{u}_k -\bar{v}_k),\tilde{u}_k -\bar{v}_k\right\rangle \bigg]\}\\
						\geq & \frac{1}{2} \mathbb{E}\Bigg\{\left\langle M\left(\tilde{y}_0-\bar{y}_0\right), \tilde{y}_0-\bar{y}_0\right\rangle
						\\&\qquad+\sum_{k=0}^{N-1}\bigg[\left\langle R_k( \tilde{u}_k-\bar{v}_k), \tilde{u}_k-\bar{v}_k\right\rangle\bigg]\Bigg\}.
					\end{aligned}
			\end{equation}}

			Due to positive definiteness of $M$ and $R_k$, we have $\tilde{u}_k=\bar{v}_k$. The uniqueness is proved.  \qed

	\section{Conclusion} \label{sec:5}

In this paper, we have investigated a class of discrete-time fully coupled nonlinear fully coupled nonlinear FBS$\bigtriangleup $E motivated by LQ  control problems. By integrating a discrete-time domination-monotonicity framework with the continuation method, we have derived essential a priori estimates and rigorously established the unique solvability of the system. To illustrate the practical relevance of our theoretical results, we have analyzed two discrete-time LQ control problems and demonstrated their solvability under the proposed framework. Future research directions include extending the model to infinite-horizon settings to explore asymptotic behaviors, formulating optimal control strategies via Pontryagin-type principles, generalizing the framework to accommodate non-Markovian or regime-switching dynamics, and developing numerical algorithms for high-dimensional systems. These advancements are expected to enhance the applicability of discrete-time stochastic control theory in financial engineering and decision-making under uncertainty.




\begin{ack}                               
This work was supported by the National
			Natural Science Foundation of China ( No.12271158), the PolyU-SDU Joint Research Center on Financial Mathematics, and the Research Centre for Quantitative Finance (1-CE03),  Hong Kong grants 15221621, 15226922 and 15225124, and partially from PolyU 4-ZZP4 and 1-ZVXA.
\end{ack}

\bibliographystyle{elsarticle-harv}

\end{document}